\documentclass[11pt,reqno]{amsart}
\usepackage{latexsym,amsthm,amssymb,amscd,amsfonts,amsmath,enumitem}
\usepackage{epsfig}
\usepackage{nicefrac}
\usepackage[all]{xy}
\usepackage{graphicx}

\usepackage[usenames]{color}

\theoremstyle{plain}
\newtheorem{thm}{Theorem}
\newtheorem*{thm*}{Theorem}
\newtheorem{theorem}{Theorem}[section]
\newtheorem{lemma}[theorem]{Lemma}                           
\newtheorem{proposition}[theorem]{Proposition}
\newtheorem{corollary}[theorem]{Corollary}
\newtheorem{cor}[thm]{Corollary}
\newtheorem*{remark*}{Remark}
\newtheorem*{remarks*}{Remarks}
\newtheorem{remark}[theorem]{Remark}
\newtheorem{remarks}[theorem]{Remarks}
\newtheorem{claim}[theorem]{Claim}

\newtheorem*{example*}{Example}
\newtheorem*{examples*}{Examples}
\newtheorem{definition}[theorem]{Definition}
\newtheorem*{definition*}{Definition}

\newtheorem*{question*}{Question}

\newtheorem*{claim*}{Claim}

\numberwithin{equation}{section}

\newcommand{\R}{\mathbb{R}}
\newcommand{\C}{\mathbb{C}}

\renewcommand{\P}{\mathbb{P}}
\newcommand{\N}{\mathbb{N}}
\newcommand{\Z}{\mathbb{Z}}

\newcommand{\T}{\mathbb{T}}

\newcommand{\CP}{\C\P}
\newcommand{\RP}{\R\P}

\newcommand{\cc}{\mathbb{\mathcal C}}
\newcommand{\cb}{\mathbb{\mathcal B}}
\newcommand{\ca}{\mathcal A}

\newcommand{\cn}{\mathcal N}

\newcommand{\cl}{\mathcal{L}}
\newcommand{\ci}{\mathcal{I}}

\newcommand{\ch}{\mathcal{H}}
\newcommand{\cp}{\mathcal{P}}
\newcommand{\co}{\mathcal{O}}
\newcommand{\cR}{\mathcal{R}}

\newcommand{\id}{\textnormal{Id}\,}

\newcommand{\nbd}{neighbourhood }
\newcommand{\nbds}{neighbourhoods }
\newcommand{\fonction}[5]
{$$ 
\begin{array}{rcccl}
 #1 & : & #2 & \longrightarrow &#3 \\
    &   & #4 & \longmapsto &#5 
\end{array}
$$}

\newcommand{\priv}{\backslash}
\newcommand{\lra}{\longrightarrow}
\newcommand{\hra}{\hookrightarrow}

\newcommand{\om}{\omega}
\newcommand{\lag}{\textnormal{Lag}}

\newcommand{\eps}{\varepsilon}
\newcommand{\sdb}{\textnormal{\small SDB}}
\newcommand{\spb}{\textnormal{\small SPB}}

\newcommand{\st}{\textnormal{st}}

\newcommand{\crit}{\text{Crit}\,}

\newcommand{\hur}{\textnormal{Hur}}
\newcommand{\wdt}[1]{\widetilde{#1}}

\newcommand{\cqfd}{\hfill $\square$\!\! \vspace{0.1cm}\\}
\newcommand{\sbull}{{\tiny $\bullet$ }}
\newcommand{\ds}{\displaystyle}
\newcommand{\im}{\textnormal{Im}\,}
\newcommand{\re}{\textnormal{Re}\,}

\newcommand{\pd}{\textnormal{\small PD}}
\newcommand{\can}{\textnormal{can}}

\newcommand{\nf}[2]{{\nicefrac{#1}{#2}}}

\newcommand{\eucl}{{\textnormal{eucl}}}
\newcommand{\cliff}{{\textnormal{cliff}}}
\renewcommand{\min}{{\textnormal{min}}}
\newcommand{\fs}{{\text{FS}}}

\newcommand{\Om}{\Omega}

\newcommand{\its}{\item[\sbull]}

\newcommand{\op}{\textnormal{Op}\,}
\newcommand{\skel}{\textnormal{Skel}}

\renewcommand{\span}{\textnormal{Span}}

\definecolor{marron}{rgb}{0.64,0.16,0.16}

\renewcommand\phi\varphi

\newcommand{\la}{\langle}
\newcommand{\ra}{\rangle}
\renewcommand{\centerdot}{{\scriptscriptstyle \bullet}}

\begin{document}
\vspace*{-2cm}
\title{Legendrian barriers in prequantization bundles}

\author{Emmanuel Opshtein}


\begin{abstract}
In \cite{opsc} we proved that the contact sphere $S^3$ has Legendrian barriers, that intersect any sufficiently long embedded  Reeb cylinders over Legendrian knots. We generalize this result to contact prequantization bundles. 
\end{abstract}

\maketitle


\section{Introduction}
Symplectic manifolds have several flavors of size, measured by so-called symplectic capacities of different natures. Dynamical capacities, like Hofer-Zehnder's, Ekeland-Hofer's, Viterbo's, or Hutchings' ECH, just to cite a few, detect  periodic orbits of Hamiltonian systems subject to various constraints and normalizations, and evaluate their symplectic actions \cite{hoze,ekho,viterbo,hutchings1}. Others measure the maximal size of natural objects that fit into the manifold. The objects used to probe the manifold can either be of the same dimension (the embedding capacities, including Gromov's width), or Lagrangian submanifolds (Cieliebak-Mohnke). 
In \cite{biran1}, Biran discovered that  removing certain Lagrangian CW-complexes from closed K\"ahler manifolds with rational symplectic classes decreases their Gromov's width, up to arbitrarily small values.  He called these Lagrangian submanifolds {\it barriers}. 
The K\"ahler assumption was soon removed by works of Donaldson and Giroux \cite{donaldson1,giroux1}. With Felix Schlenk,  we proposed an extension of Biran's framework to the open case,  and studied it  in the simplest example of the $4$-dimensional bidisc \cite{opsc}. The result was that in this context, the removal of some explicit isotropic CW-complexes decreases {\it all} symplectic capacities. 

By comparison, the quantitative aspects of contact geometry are more subtle. For instance, the group of contactomorphisms allows for arbitrary contractions of a Darboux ball, but there still exists a contact non-squeezing phenomenon \cite{elkipo,sandon1,frsazh,cant}. In parallel with our work on the rigidity of Lagrangian skeleta, a  Legendrian barrier phenomenon in the contact $3$-sphere was shown: 
\begin{thm*}[\cite{opsc}]
Let $\alpha_0$ be the contact form on $S^3$ obtained by restricting the form $R_1d\theta_1+R_2d\theta_2$ to $S^3(1)\subset \R^4$. There exists a Legendrian CW-complex $\Lambda_k$ such that, for any contact form $\alpha\leq \alpha_0$, any closed Legendrian submanifold $\Lambda\subset S^3$ whose $\alpha$-Reeb flow remains in $S^3\priv \Lambda_k$ for a time $\frac 1k$ has an $\alpha$-Reeb chord of length $\leq \frac 1k$.
\end{thm*}
The aim of the present paper is to generalize this statement to  contact prequantization bundles. Before we state our main result, we need to introduce the notation and some concepts. Let $(N,\tau)$ be a symplectic manifold with integral symplectic class (so $[\tau]\in H_2(N,\Z)$). There are several notions of polarizations. The most basic one, a pair  $(\Sigma,\lambda)$ consisting of a smooth symplectic hypersurface Poincaré-Dual to $k[\tau]$ in $N$ and a Liouville form on $N\priv \Sigma$ is too loose for our purpose. Throughout this paper, a  polarization of 
$(N,\tau)$ is a richer data obtained as follows. Consider first a complex line bundle $(\cl,\nabla)\to N$ with Chern-class $-[\tau]$ equipped with a Hermitian connection, as well as its tensor products $\cl^k$. Once a compatible almost complex structure $J$ is chosen on $(N,\tau)$, one may look for sections that behave nicely with respect to the complex structures on the base and the fiber, respectively. When $(N,\tau,J)$ is K\"ahler for instance, if $\cl$ is moreover  holomorphic, Kodaira's theory guarantees the existence of many holomorphic sections at least for large $k$, some vanishing transversally. Kodaira's theory has been extended to the symplectic category by Donaldson and Giroux in the following form: when $k$ is large, there exist quasi-holomorphic sections $s_k$ of $\cl^{k}$ that vanish transversally. Their vanishing loci are smooth symplectic hypersurfaces of $(N,\tau)$ whose complements are naturally endowed with Liouville forms which are gradient-like for $-\ln|s_k|$ \cite{donaldson1,giroux1,giroux2}. This procedure (whether in the K\"ahler setting, or in Donaldson-Giroux's) is in fact the only known source of smooth symplectic hypersurfaces with Weinstein complements. Although we still denote them $(\Sigma,\lambda)$, our polarizations must be understood in the previous sense: $\Sigma$ is the transverse vanishing set of a holomorphic or a quasi-holomorphic section $s$ of $\cl^k$, for some $k$ called the degree of the polarization, and $\lambda$ is the corresponding Liouville form on $N\priv \Sigma$. Such a polarization comes with an isotropic {\it skeleton}, noted $\skel(\Sigma,\lambda)$,  defined as the maximal compact invariant subset of $X_\lambda$ in $N\priv \Sigma$. Skeleta are known to be critical isotropic when $k\geq 2$ \cite{gespze}. When $\ln| s |$ is Morse on $N\priv \Sigma$, this skeleton is the CW-complex formed by all the stable manifolds of the critical points of $-\ln|s|$ under the flow of $\lambda$. Biran's Lagrangian barriers are such skeleta.


We now define the {\it contact prequantization bundle} associated to an integral symplectic manifold $(N,\tau)$. Let $\cp\to N$ be the $S^1$-bundle over $N$ with Euler-class $-[\tau]$ and $\alpha$ a connection form for this bundle with curvature $\tau$. This means that 
$$
\left\{\begin{array}{l} \alpha_{|S^1}=d\theta\\ d\alpha=\pi^*\tau.\end{array}\right.
$$ 
Although $\alpha$ is not completely determined by the above conditions, different choices of connections lead to strictly contactomorphic forms. The contact prequantization associated to $(N,\tau)$ is the contact manifold $\cp(N,\tau):=(\cp,\alpha)$, and we insist that it is endowed with a contact form, not just a contact structure. 

Recall finally  that in a contact manifold $(M,\xi=\ker \alpha)$, any function 
$H\in \cc^\infty(M\times[0,1])$ gives rise to a contact Hamiltonian vector field, defined by 
$$
\left\{\begin{array}{l}
X_H=HR_\alpha+Z,\hspace{,3cm} Z\in \xi,\\
d\alpha(Z,\cdot)=-dH_{|\xi}. 
\end{array}\right.
$$
The map $\Phi_H^t:=\Phi_{X_{H}}^{[0,t]}$ is a contact diffeomorphism, meaning that it preserves the contact distribution $\xi$. 
Given a time-dependent function $H$ and two subsets $X,Y\subset M$, we say that $H$ has a chord of length $\ell$ between $X$ and $Y$ if there exists a  point  $x\in X$ such that $\Phi^\ell_H(x)\in Y$. Notice that by convention, the trajectory is required to start on $X$ at time $0$. Such a trajectory is called an $H$-chord from $X$ to $Y$.

\begin{thm}\label{t:legbarp} Let $(N^{2n},\tau)$ be a closed symplectic manifold with integral symplectic class and $(\cp^{2n+1},\alpha)\overset{\pi}\to (N,\tau)$ be its prequantization bundle. Let $\Gamma_k\subset N$ be the isotropic skeleton associated  to a quasi-holomorphic polarization of degree $k\geq 2$ of $(N,\tau)$. Let $\Lambda_k$ be a connected component of a Legendrian lift of $\Gamma_k$ to $\cp$.  Then, 
\begin{itemize}
\its $\Lambda_k$  is a compact subset of $\cp$, that intersects each circle $\pi^{-1}(p)$, $p\in \Gamma_k$ in exactly $k$ (equi-distributed) points,
\its For any closed Legendrian submanifold $\Lambda\subset (\cp,\alpha)$  and any smooth contact Hamiltonian $H:\cp\times [0,1]\to [1,+\infty)$, 
either there exists an $H$-chord from $\Lambda$ to $\Lambda_k$ of length $\leq \frac 1k$, or the set 
$$
\bigcup_{t\in [0,\frac 1k]} \Phi^t_H(\Lambda)
$$
is not embedded. 
\end{itemize} 
When $H$ is autonomous, the second bullet can be stated alternatively:
\begin{itemize}
\its There exists an $H$-chord from $\Lambda$ to $\Lambda\cup \Lambda_k$ of length $\leq \frac 1k$. 
\end{itemize}
\end{thm}
\noindent In other terms, like several results on the symplectic rigidity of the lifts of natural rational Lagrangian manifolds, {\it e.g.} \cite{disu2,kilgore}, we get:
\begin{center}
{\it The Legendrian lift of a Lagrangian barrier is a Legendrian barrier.} 
\end{center}
Indeed, its removal obstructs the existence of certain objects that exist in $\cp(N,\tau)$. For instance:
\begin{cor}\label{c:legbar}
There exists no contact form $\alpha'\leq \alpha$ on $(\cp(N,\tau),\xi=\ker \alpha)$ with a strict contact embedding 
$$
(\overline{D(\nf 1k)}^n\times S^1,\lambda_\st+d\theta)\hra (\cp(N,\tau)\priv \Lambda_k,\alpha'). 
$$
\end{cor}

The first bullet of theorem \ref{t:legbarp} relies on the fact that, although not exact,   Lagrangian skeleta are rational: the areas of the discs they bound are multiples of $\nf 1k$. This fact is probably well-known to the experts, but we provide a precise statement and proof in proposition \ref{p:exskel}. Let us now comment on the main content of theorem \ref{t:legbarp}, its second bullet.

In \cite{enpo}, Entov and Polterovich introduced a notion of interlinkedness. Two Legendrian submanifolds $\Lambda_1,\Lambda_2\subset (M,\xi=\ker \alpha)$ are $\delta$-interlinked for $\alpha$ if all contact Hamiltonian vector fields associated to smooth functions $\geq 1$  have a chord from $\Lambda_1$ to $\Lambda_2$ of length at most $\delta$. Although interlinkedness does not depend on the contact form $\alpha$, the constant of interlinkedness $\delta$ does: expanding $\alpha$ by a factor $C$ shrinks $\delta$ by the same factor. This notion has been studied for instance in \cite{alar,disu}. In some sense, the Legendrian barriers $\Lambda_k$ of theorem \ref{t:legbarp}  are {\it universal interlinkers}.  Let us dedicate some lines to propose a proper definition of this notion. We would like to say that a subset $X\subset (M,\xi=\ker \alpha)$ is a universal $\delta$-interlinker if $X$ is $\delta$-interlinked with every Legendrian submanifold $\Lambda$. Thinking a bit shows the lack of relevance of such a definition. Since there exist local Legendrian knots in arbitrary \nbds of points of $M$, a universal interlinker $X$ would have to intersect {\it all} Reeb orbits (associated to the Hamiltonian $H\equiv 1$), hence would need to be of codimension $1$, and the property would in fact have no relation with the Legendrian nature of $\Lambda$.  Still, if we force the test Legendrian submanifolds to be non-local, or more precisely not localized near a Hamiltonian trajectory, the interlinking property has a chance to hold. As in \cite{opsc}, we characterize the non-locality of a Legendrian submanifold by requiring that it has no short self-chord. 

\begin{definition*} A subset $X$ of a contact manifold $(M,\xi=\ker \alpha)$ is a $(\delta,\delta')$-universal interlinker with respect to $\alpha$ if for any Legendrian submanifold $\Lambda\subset M$ and  smooth  function $H\in \cc^\infty(M\times[0,1])$ with $H\geq 1$, 
\begin{itemize}
\its either there exists an $H$-chord from $\Lambda$ to $X$ of length less than $\delta$,
\its or there exists $t_0\in [0,\delta']$ and an $H$-chord from $\Lambda_{t_0}:=\Phi_H^{t_0}(\Lambda)$ to itself starting at $t_0$ and getting back to $\Lambda_{t_0}$ before  time $\delta'$. 
\end{itemize}
The second bullet can be restated by saying that the $\delta'$-sweep of $\Lambda$ along $H$, defined by 
$$
\bigcup_{t\in [0,\delta']}\Phi^t_H(\Lambda)
$$
is not embedded. When $H$ is autonomous, it even simplifies to:  
\begin{itemize}
\its there exists an $H$-chord from $\Lambda$ to itself of length $\leq \delta'$. 
\end{itemize}
\end{definition*}
As above, the {\it notion} of universal interlinkedness depends only on the contact structure but the exact values $(\delta,\delta')$ depend on the contact form.  Testing this definition with autonomous Hamiltonians is already interesting. In this case, the Hamiltonian flow is exactly the Reeb flow of the form $\alpha':=H^{-1}\alpha$, the function $H^{-1}$ taking values in $(0,1]$. Moreover, the class $\cl(\alpha',\delta')$ of Legendrian submanifolds that have no self $\alpha'$-chord of length less than $\delta'$ is invariant under the $\alpha'$-Reeb flow. Thus, $X$ being a $(\delta,\delta')$-universal interlinker for the form $\alpha$ guarantees in particular the following: for any contact form $\alpha':=f\alpha$ with $0<f\leq 1$, for any $\Lambda \in\cl(\alpha',\delta')$, there exists an $\alpha'$-Reeb chord from $\Lambda$ to $X$ of length $\leq \delta$.
With this terminology, the content of theorem \ref{t:legbarp} is that $\Lambda_k$ is a $(\frac 1k,\frac 1k)$-universal interlinker. It might be interesting to understand the dependency between the constants $\delta$, $\delta'$ of the definition: does there exist $(\delta,\delta')$-universal interlinkers with sharp constants $\delta\neq \delta'$ ?


Some hard tools based on generating functions or Legendrian Contact Homology may detect Reeb chords.  The latter has even been computed for examples of Legendrian lifts in prequantizations \cite{digo}. The  non-smoothness of $\Lambda_k$ puts it however  outside the usual scope of these techniques. Indeed, although Legendrian lifts of Lagrangian immersions may well be embedded, the lifts of our CW-complexes are not. This already happens for the Legendrian lift of the skeleton of $S^2$ consisting of half great circles that cut $S^2$ into $k$ equal discs \cite{opsc}.  
Our proof consists of an embedding construction à la Mohnke that reduces the statement to a problem on estimating the Lagrangian area class of some Lagrangian embedding associated to $\Lambda$. As in Biran's initial work, the embedding takes place in a space that somehow forgets about $\Lambda_k$, so that the hard part of the proof avoids the investigation of the latter. Still, recent results on detecting rigidity properties of  non-embedded isotropic skeleta {\it via} symplectic homology tools \cite{tova} give hope that these Legendrian barriers can be studied by other means than the one developped in this paper, leading to further rigidity properties. For instance, our theorem  gives a bound on a "contact Gromov width" of $\cp(N,\tau)\priv \Lambda_k$ (corollary \ref{c:legbar}).  But as for Lagrangian skeleta,  the complements of the $\Lambda_k$ may get small whatever  contact invariant is used for measure. Let us ask a precise question in this direction. In \cite{cosa}, V. Colin and S. Sandon introduce the Legendrian discriminant length between isotopic Legendrian submanifolds, defined as follows. Say first that a Legendrian isotopy $(\Lambda_t)_{t\in [0,1]}\subset (M,\xi)$ is embedded if 
$$
\bigcup_{t\in [0,1]}\Lambda_t
$$
is an embedding of $[0,1]\times \Lambda$ into $M$. Then define the discriminant length of a Legendrian isotopy as the minimal number $n$ for which $(\Lambda_t)_{t\in [0,1]}$ may be written as a concatenation of $n$ embedded isotopies. And finally define the  discriminant distance between two isotopic Legendrian submanifolds as the shortest length of a Legendrian isotopy between them. Theorem \ref{t:legbarp} shows that any isotopy of a Legendrian $\Lambda$ that remains in the complement of $\Lambda_k$ through the Hamiltonian flow of a function larger than $1$ has length at least $k$. Hence the following question:
\begin{question*}
Does there exist a sequence of Legendrian pairs $(\Lambda_1^k,\Lambda_2^k)$ in $\cp\priv \Lambda_k$ whose discriminant distance in $\cp\priv \Lambda_k$ grows much faster than in $\cp$ ? 
\end{question*}

Let us finally state a rigidity result on Lagrangian skeleta:
\begin{thm}\label{t:rigskel} Let $(M,\om)$ be a symplectic manifold with integral symplectic class and 
$\Gamma_k\subset M$ be the skeleton of a polarization of degree $k$. Then any closed Lagrangian submanifold $L$ in $M\priv \Gamma_k$  bounds a symplectic disc in $M\priv \Gamma_k$ of area $<\nf 1k$. In particular, for the Cieliebak-Mohnke capacity we have: 
$$
c_\lag(M\priv \Gamma_k)\leq \frac 1k. 
$$ 

\end{thm}
The proof  is rather straightforward after Biran's original paper and Cie\-liebak-Mohnke's \cite{biran2,cimo}. We take advantage of the present paper to provide a full proof, the road for theorem \ref{t:legbarp} passing very nearby. Interestingly, theorem \ref{t:legbarp} follows from theorem \ref{t:rigskel}, but then implies an extension of the latter for some Lagrangian immersions: 
\begin{cor}\label{c:rigskel}
Let $(M,\om)$ be a simply connected symplectic manifold with integral symplectic class and 
$\Gamma_k\subset M$ be the skeleton of a polarization of degree $k$. Let $L$ be the image of a Lagrangian immersion in $M\priv \Gamma_k$ that has only transverse self-intersections  and verifies $[\om]\cdot H_2(M,L)\subset \frac 1N\Z$. Then $N\geq k$. 
\end{cor} \vspace{,2cm}

\paragraph{\bf Organization of the paper.}
In section \ref{s:examples}, we illustrate theorems \ref{t:legbarp} and \ref{t:rigskel} by  specific examples where the skeleta can be computed, and we compare their statements with existing results. In section \ref{s:idea}, we explain the proof of theorem \ref{t:legbarp} in a very particular case, deduce corollary \ref{c:rigskel}, and try to give the general idea of the actual proof. In section \ref{s:tech} we gather several general statements that may be found scattered throughout the litterature. For the sake of readability  and to keep the exposition  self-contained, we state these results in the  most suitable form for our proof of theorem \ref{t:legbarp} and provide at least brief proofs. We also take a short detour to prove theorem \ref{t:rigskel}. In section \ref{s:sdb} we set the framework of our proof, discussing symplectic disc and projective bundles and their relations with prequantizations. This section is mainly meant for setting the notation, save for the more original paragraph \ref{s:spbb}. The sections \ref{s:sec}, \ref{s:lvf} and \ref{s:taming} are the technical (and mostly computational) core of the paper. \vspace{,2cm}

\paragraph{\bf Aknowledgements.} I wish to thank Felix Schlenk for the very nice collaboration \cite{opsc}, in which he had this inspiring idea of Legendrian barrier phenomenon. This papers owes a lot to him, through this first collaboration,  his enthusiasm for the generalization proposed here, an his careful re-reading of the first version of this paper. I also thank warmly Margherita Sandon for the time she took to answer my several questions on contact geometry, and for the interest she has shown in this work. \vspace{,2cm}

\paragraph{\bf Notation and convention.}

\begin{itemize}
\its The image of the Hurewicz morphism $\pi_2(M)\to H_2(M,\Z)$ is denoted $H_2^\hur(M,\Z)$. 
\its Given a subset $S$ of a manifold $M$,  $\op(S,M)$ stands for an arbitrarily small \nbd of $S$ in $M$.  We sometimes write $\op(S)$ if $M$ is understood from the context. We sometimes need pointed neighbourhoods, that we denote $\op(S)^*:=\op(S)\priv S$.
\its For a symplectic manifold $(M,\om)$, which may be a symplectic submanifold of an ambient one, we write $\ca_\om:H_2(M,\Z)\to \R$ for the evaluation of $[\om]$. 
\its Similarly, when $L\subset (M,\om)$ is a Lagrangian submanifold, we call $\ca_\om^L:H_2(M,L,\Z)\to \R$ its area morphism. When $\im \ca_\om^L$ is a discrete subgroup of $(\R,+)$, we call its positive generator the area class of $L$ and denote it $\ca_\min(L)$. 
\its Our norms are always written with simple bars $|\cdot|$. No vector space is endowed with more than one norm, so we never distinguish our  norms in the notation.
\its The Fubini-study form $\om_\fs$ on $\CP^n$ is normalized by giving area $1$ to projective lines. 
\end{itemize}

\section{Examples}\label{s:examples}
It is natural to ask how theorems \ref{t:legbarp} and \ref{t:rigskel} fit in the landscapes of non-removable intersections in symplectic geometry  and of existence of Legendrian chords in contact geometry. In high degrees, theorem \ref{t:rigskel} provides non-removable intersections at small scales, which hardly compares to anything already investigated in symplectic geometry except in \cite{posh}. At the monotone scale however, many conjectures and results state the existence of non-removable intersections between a large class of Lagrangian submanifolds and a fixed object. For instance, it is conjectured that all monotone Lagrangian tori (or even maybe submanifolds) $L\subset \CP^n$ intersect the Clifford torus, and this conjecture is known to hold for any heavy Lagrangian submanifold \cite{enpo2}.  Heaviness is in turn achieved for all monotone $L$ whose (well-defined) Floer homology with {\it some} coefficients do not vanish. 
In the contact side, despite  similar wording, the phenomenon of universal interlinking
does not really compares to Entov-Polterovich's work, in the sense that they are concerned with finding conditions on pairs of objects that guarantee their interlinking, while theorem \ref{t:legbarp} produces one object which is interlinked with a large class of Legendrian submanifolds. Forgetting its quantitative nature, theorem \ref{t:legbarp} is reminiscent of Givental's result of existence of Legendrian chords between $\R\P^n$ and its Legendrian isotopies in $\RP^{2n+1}$  \cite{givental}. The aim of this section is to situate the present paper in this broader context, so we state examples of what can be derived from the present work in these directions. We do not try to provide proofs or justifications of all our statements.  In particular the lengthy computations of skeleta are skipped.

Before getting to the point, let us finally  add a word of warning. The only example of smooth skeleton associated to a smooth polarization of degree $\geq 2$ known to the author is   $\R\P^n\subset \C\P^n$, which is the skeleton of a polarization of degree $2$.  The examples below therefore claim non-removable intersections between two subsets, one of which is usually non-smooth, albeit a CW-complex unless explicitely stated. \vspace{,3cm}

\paragraph{\bf Illustration of theorem \ref{t:legbarp}.}
Let us  start with the polarization of $(\CP^n,\om_\fs)$ of degree $2$ whose skeleton is $\RP^n$. 
 The prequantization of $(\C\P^n,\om_\fs)$ is  $S^{2n+1}(1)$, and the Legendrian lift of $\R\P^n$ is $S^{2n+1}_\R:=\R^{n+1}\cap S^{2n+1}(1)\subset \C^{n+1}$. A contact form $\alpha\leq \alpha_\st$ corresponds to the contact form induced by the standard Liouville form on $\C^n$ on the boundary of some star-shape domain $\Om\subset B^{2n+2}(1)$, $S^{2n+1}_\R$ corresponding to $\partial \Om_\R:=\partial \Om\cap \R^{n+1}$. Theorem \ref{t:legbarp} guarantees that any closed Legendrian submanifold of $\partial \Om$ either has a self-Reeb chord  or a Reeb-chord to $\partial \Om_\R$, of length $<\nf 12$. Forgetting about the length estimate, this is a weaker result than Mohnke's, which asserts that {\it every} Legendrian submanifold has a self Reeb-chord \cite{mohnke}. But it shows that a class of Legendrian submanifolds, those that have no too short Reeb-chords, also have a Reeb chord to $S^{2n+1}_\R$, which is reminiscent of Givental's result cited above. 

The prequantization of $(\C\P^n,2\om_{\fs})$ is the contact manifold $\R\P^{2n+1}$, and the Legendrian lift of $\RP^n$ is $\R\P^n$. It would be tempting  to apply theorem \ref{t:legbarp} in order to get a result à la Givental. This is however not possible because $\R\P^n\subset (\C\P^n,2\om_\fs)$ is the skeleton of a polarization of degree $1$, for which theorem \ref{t:legbarp} does not apply. A skeleton of a polarization of degree $2$ in $(\C\P^n,2\om_\fs)$ corresponds to the skeleton of a polarization of degree $4$ in $(\C\P^n,\om_\fs)$. In \cite{biran2}, Biran computes examples of such skeleta:
$$
\Gamma_2(\C\P^n,2\om_\fs)=G\cdot \R\P^n,
$$
where $G:=(\Z_2)^n\subset PU(n)$. The orbit $G\cdot \R\P^n$ consists of a specific collection of $2^{n-1}$ real projective planes. Our theorem then implies that for any contact form $\alpha\leq \alpha_\st$ on $\R\P^{2n+1}$,  a Legendrian submanifold either has a self Reeb-chord, or a Reeb chord to $G\cdot \R\P^n$ of length $<\nf 12$. Again, restricting the class of Legendrian submanifolds to the ones that are Legendrian isotopic to $\R\P^n$ and forgetting about the length estimate gives a poor result compared to Givental's, that this Legendrian submanifold  has a Reeb chord to {\it one} $\R\P^n$. 
\vspace{,3cm}

\paragraph{\bf Illustration of theorem \ref{t:rigskel}.}
We start again with the pair $(\CP^n,\om_\fs)$ with skeleton $\RP^n$ in degree $2$. Then theorem \ref{t:rigskel} tells that any Lagrangian submanifold in $\C\P^n\priv \R\P^n$ has an attached disc of area $<\nf 12$. A small refinement of the argument provided for theorem \ref{t:rigskel} even proves that the same holds for any Lagrangian submanifold in $B^{2n}(1)\priv \R^n$  (see theorem \ref{t:rigskelaff}). So for instance a torus $T(x)$ cannot be squeezed in $B^{2n}(1)\priv \R^n$ when $x\geq \nf 12$. Since however such a torus cannot  be squeezed in $B^{2n}(1)$ already when $x\geq \nf 1n$ \cite{cimo}, theorem \ref{t:rigskel} is void in this case.

In $(\C\P^n,\om_\fs)$, the Clifford torus is the skeleton of a {\it singular} polarization of degree $n+1$, associated to the section $s$ of $\co(n+1)\to \P^n$ corresponding to the homogeneous polynomial $z_0z_1\dots z_n$. If we regularize this section by considering instead the section $s_\eps$ associated to the polynomial $z_0\dots z_n+\eps(z_0^{n+1}+\dots +z_n^{n+1})$, one checks that for $\eps\ll 1$,  $s_\eps$ has exactly $n+1$ critical points obtained from
$p_0:=[1:0:\dots:0]$ by permutations of coordinates, besides the critical points of $s$. These new critical points are not Morse, but have computable stable sets, defined as the set of points attracted by these critical points under the flow of $-\nabla \ln |s_\eps|$: 
$$
W^s(p_0)=\{R_1=\dots=R_n<\nf 1{n+1},\; \theta_1+\dots+\theta_n=0\}\subset \C^n\approx\P^n\priv\{z_0=0\}. 
$$ 
This is a Lagrangian cone over an $n-1$-dimensional torus with boundary attached to the Clifford torus 
$$
\T_\cliff=\{R_1=\frac 1{n+1},\dots,R_n=\frac 1{n+1},\theta_1,\dots,\theta_n\}
$$
along the linear sub-torus $\theta_1+\dots+\theta_n=0$. Except in dimension $4$, when $s_\eps$ is in fact Morse, this stable set has a cone-type singularity at $p_0$, and is not a CW-complex. Taking into account the symetries of our section, denoting $u_i$ the permutation of coordinates between $z_0$ and $z_i$, and putting $p_i:=u_i(p)=\{z_j=0\;\; \forall j\neq i\}$, we get
$$
C_i:=W^s(p_i)=u_i\big( W^s(p_0)\big). 
$$
The full skeleton of the polarization by $s_\eps$ is therefore
$$
\Gamma_3(\C\P^n)=\T_{\cliff}\cup \big(\bigcup_{i=0}^{n} C_i\big).  
$$
Theorem \ref{t:rigskel} then tells that a monotone Lagrangian submanifold of $\C\P^n$ either intersects the Clifford torus, or at least one of the $C_i$. This is obviously  less than the conjecture that any monotone Lagrangian submanifold has an intersection with the Clifford torus, which holds for heavy Lagrangian submanifolds \cite{enpo2}. 

\begin{figure}[h!]
\begin{center} 
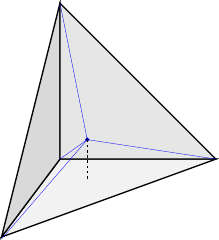
\caption{Toric projection of the skeleton $\Gamma_3(\CP^n)$ (in blue).}
\label{fig:skeleton}
\vspace{,1cm}
\begin{minipage}{0,8\textwidth}{\footnotesize $\Gamma_3(\CP^n)$ is the union of the Clifford torus and  Lagrangian lifts of the blue lines. Above each point of these lines, the lift is a linear sub-torus of the fiber.}\end{minipage}
\end{center}
\end{figure}

The same description can be carried on for the skeleton of a polarization of degree $2$ of $S^2(1)^n$, obtained by a symetric perturbation of a section whose zero set is the union of the hypersurfaces $\{z_i\in\{0,\infty\}\}$. The  skeleton is the union of the Clifford torus and $2^n$ Lagrangian cones over $n-1$-dimensional tori that all abute to the Clifford torus. Namely, 
$$
\Gamma_2(S^2(1)^n)=\T_\cliff\cup(\bigcup_{i=1}^{2^n} C_i), 
$$
where 
$$
C_0=\{(R,\dots,R,\theta_1,\dots,\theta_n)\;|\; R\in [0,\nf 12],\; \theta_1+\dots+\theta_n=0\}
$$
in toric coordinates, and the $C_i$ are obtained from $C$ by permutations of the factors and inversions in each factor.  
Theorem \ref{t:rigskel} states that $\Gamma_2(S^2(1)^n)$ intersects any monotone Lagrangian submanifold. In this case again, this  is  weaker than the natural conjecture, which holds for heavy Lagrangian submanifolds  \cite{enpo2} and for monotone tori in dimension $4$ \cite{hike}. 

\begin{figure}[h!]
\begin{center} 
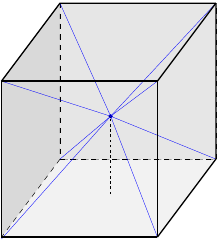
\caption{Toric projection of the skeleton $\Gamma_2(S^2(1)^n)$ (in blue).}
\label{fig:skeleton2}
\end{center}
\end{figure}

Finally, theorem \ref{t:rigskelaff} gives similar results in  the balls or the polydiscs. For instance, consider the set $\Gamma_3(\CP^n)$ defined right above, and transform it to a set $\Gamma_3(B^{2n}(1))\subset B^{2n}(1)$ by first applying an element $u\in \P U(n)$ that takes the hyperplane $\ch:=\{z_0+\dots+z_n=0\}$ to $\ch_\infty:=\{z_0=0\}$ and  then removing $\ch_\infty$. Then the tori $\T^n(x,\dots,x)$, with $x\geq \frac 1n$ cannot be displaced from $\Gamma_3(B^{2n}(1))$ in $B^{2n}(1)$. When $\nf 1{n+1}<x<\nf 1n$, these tori can be taken into the complement of $u(\T_\cliff)$ but have to intersect one of  the $u(C_i)$. A maybe less complicated way of saying the same is the following: in $\C\P^n\priv \ch$, there are Lagrangian tori whose attached discs that do not intersect $\ch$ all have area a multiple of $x$, for $x\in (\nf 1{n+1},\nf 1n)$. These Lagrangian tori cannot be made disjoint from $\Gamma_3(\C\P^n)$ in $\C\P^n\priv \ch$.

\section{Idea of proof}\label{s:idea}
The case $H=1$ of theorem \ref{t:legbarp} is already non-trivial but can be explained easily, up to a certain point at least, by an argument  very similar to the one  in \cite[Section 8]{opsc}. In this case, the Hamiltonian vector field is the Reeb vector field on $(\cp,\alpha)$, which generates the $S^1$-action of the circle bundle $\cp\overset\pi\lra N$. Saying that a Legendrian $\Lambda$ has no Reeb chord to $\Lambda_k$ of length $\leq \nf 1k$ amounts to saying that it does not intersect the $-\nf 1k$-sweep of $\Lambda_k$ under the Reeb flow, defined by 
$$
\bigcup_{t\in[-\nf 1k,0]} \Phi_{R_\alpha}^t (\Lambda_k). 
$$ 
Since $\Lambda_k$ consists of  $k$ equi-distributed points in each fiber circle over the points of $\Gamma_k$, this $-\nf 1k$-sweep is exactly $\pi^{-1}(\Gamma_k)$, so the absence of Reeb chord from $\Lambda$ to $\Lambda_k$ of length $\leq \nf 1k$ means that $L:=\pi(\Lambda)$ lies in $N\priv \Gamma_k$. In dimension $3$, $N$ is a surface that $\Gamma_k$ cuts into discs of area $\nf 1k$ (see \cite[Section 3.1]{opsc} or proposition \ref{p:exskel}). Thus $L$ is a Lagrangian immersion into a $2$-disc of area $\nf 1k$. It therefore bounds a disc of area $<\nf 1k$, and any Legendrian lift of the boundary of this disc starting at a point in $\Lambda$ provides a path whose endpoints are on the same fiber, thus  a Reeb self-chord of $\Lambda$. The   length of this Reeb-chord is the area  of the disc, hence the result. In higher dimension when $L$ is embedded, theorem \ref{t:rigskel} provides a disc of area $< \nf 1k$ bounded by $L$ that can be used exactly as in dimension $3$ to produce the short Reeb chord we are seeking. When $L$ is only immersed, this easy argument does not seem to go through, and it seems that the long way is required. Of course, when $H$ is not constant, the Reeb flow is not periodic any more, and the above argument breaks down much before this point. This discussion also explains corollary \ref{c:rigskel}. Indeed, if $L\subset (M,\om)$ is a Lagrangian immersion with transverse self-intersections in a simply connected integral symplectic manifold with  $\ca_\min(L)\in \frac 1N\Z$, a minimal Legendrian lift of $L$ to $\cp(M,\om)$ is a Legendrian embedding $\Lambda \subset \cp(M,\om)\priv \pi^{-1}(\Gamma_k)$ whose intersections with a given fiber lies in an equi-distributed subset of $N$ points. Thus $\Lambda$ has no $\alpha$-Reeb chord to $\Gamma_k$, and no $\alpha$-self Reeb chord of length $\leq \nf 1N$. Applying theorem \ref{t:legbarp} to $\Lambda$, $\Gamma_k$ and $H=1$ shows that $\nf 1N<\nf 1k$, so $N>k$. 

Let us now describe briefly the actual proof of theorem \ref{t:legbarp}. In the general framework of a time dependent contact Hamiltonian $H$ on $(\cp,\alpha)=\cp(N^{2n},\tau)$, instead of projecting the picture down from $(\cp,\alpha)$ to $(N,\tau)$, we rather lift it to the symplectization $S\cp(N,\tau):=(\cp\times \R_*^+,R\alpha)$ and use Mohnke's trick. The latter associates to an embedded $T$-sweep of $\Lambda$ in $\cp\priv \Lambda_k$ a Lagrangian submanifold $L$ in $S\cp(N,\tau)\priv (\Lambda_k\times \R)$ with $\ca_\min(L)=T$. The assumption $H\geq 1$ guarantees  that $L\subset S^-\cp(N,\tau):=\cp\times (0,1)$. The latter  is symplectomorphic to the punctured symplectic disc bundle $\sdb(N,\tau)\priv N\to(N,\tau)$. The technical core of the paper consists in computing a polarization $(\Sigma_k,\lambda_k)$ of degree $k$ on $\sdb(N,\tau)$ whose skeleton is precisely $\Lambda_k\times [0,1)$. Since $L$ does not intersect this skeleton, Biran's decomposition theorem  shows that it lies in the symplectic disc bundle over $\Sigma_k$ of radius $\nf 1k$, which compactifies to a ruled symplectic manifold of width $<\nf 1k$. The Lagrangian capacity of the latter can be shown to be bounded by $\nf 1k$ (this is how theorem \ref{t:rigskel} is proved). Thus, 
$$
T=\ca_\min(L) <\nf 1k,
$$
so no $H$-sweep of length larger than $\nf 1k$ can be embedded into $\cp\priv \Lambda_k$.

Beside the main difficulty of uncovering a polarization with computable skeleton, the secondary difficulty in this line of argument is to get a correct notion of polarization of $\sdb(N,\tau)$, which is an open symplectic manifold. The concept of Liouville polarization introduced in \cite{opsc} is not  relevant here. Our approach consists in compactifying $\sdb(N,\tau)$ to a sphere bundle, which is unfortunately not symplectic. Biran's decomposition theorem must then be adapted to this context, and the analysis requires $k\geq 2$. This last assumption is not of technical nature. The result does not hold for $(N,\tau)=(\P^n,\om_\fs(1))$ and $k=1$, because then $\Gamma_1$ is a point, and $\Lambda_1$ one Hopf circle. 

\section{Technical preliminaries}\label{s:tech}
We collect in this section a certain number of technical results, that may be considered folklore. 

\subsection{Lagrangian capacity of ruled symplectic manifolds.} 
Theorem \ref{t:legbarp} is based on the following rigidity result. We recall that $H_2^\hur(Y)$ denotes the image of the Hurewitz morphism $\pi_2(Y)\to H_2(Y,\Z)$.  
\begin{theorem}\label{t:cimo}
Let $(S^2,\om_{\fs}(a))\hra (M^{2n},\om)\overset\pi\lra Y^{2n-2}$ be a closed symplectic manifold with the structure of a symplectic sphere bundle, the parameter $a\in \R$ being the area of the fibers. Assume the existence of two  disjoint sections $Y\overset {\sigma_{0,\infty}}\lra \Sigma_{0,\infty}\subset (M^{2n},\om)$  such  that $\sigma_\infty^*\om$ is a symplectic form, and denote $X:=M\priv \Sigma_\infty$. If 
$$
a<\im(\ca_{\sigma_0^*\om}:H^\hur_2(Y)\to \R)\cap \R^*_+,
$$
then any closed Lagrangian sumanifold $L\subset (X,\om)$ 
has an attached symplectic disc 
$$
u:(D,\partial D)\lra (X,L)
$$
of area $\leq a$. 
\end{theorem}
 Theorem \ref{t:legbarp} follows from the usual SFT-type argument, as used for instance in \cite{cimo}, and is based on the following even more classical lemma:
\begin{lemma}\label{l:holfib} 
Under the assumptions of theorem \ref{t:cimo}, let $V'\Subset V$ be two \nbds of $\Sigma_\infty$ and $p\in M\priv V$.  For any $\om$-compatible almost complex structure $J$ on $M$ for which $\Sigma_\infty$ is $J$-holomorphic, there exists an arbitrarily close perturbation $J'$ of $J$ in the $\cc^1$-topology that coincides with $J$ on $(M\priv V)\cup V'$ and a pseudo-holomorphic sphere 
$$
u:(S^2,j)\lra (M,J')
$$ 
that verifies $[u]=[\pi^{-1}(\{\star\})]=:F\in H_2(M)$ and $p\in \im(u)$. 
\end{lemma}
This lemma directly follows from Gromov's theory of pseudo-holomorphic curves \cite{gromov,mcsa2}. A full proof  under slightly different assumptions can be found for instance in \cite[Theorem 5.3.A]{biran1}. The difference with our context is a trading between Biran's assumptions {\bf (L)} and ours on the existence of the symplectic sections $\sigma_0,\sigma_\infty$ and the area estimate of the fibers. Both aim at forcing the compactness of the moduli space of $J$-holomorphic spheres in the class of the fiber, that we call $F$. In our context, if $u_n:S^2\to M$ is a sequence of $J$-holomorphic spheres with $[u_n]=F$, Gromov's compactness theorem gives  a subsequence that converges in Gromov's sense to a bunch of $J$-holomorphic spheres $v_1,\dots,v_k:S^2\to M$ with $\sum [v_i]=F$.  Decomposing 
$$
H_2(M,\R)=\la F\ra+ \sigma_0(H_2(Y,\R)),
$$
write $[v_i]=\alpha_iF+B_i$. Since $\Sigma_0\cap \Sigma_\infty=\emptyset$, $B_i\cdot [\Sigma_\infty]=0$, so $[v_i]\cdot [\Sigma_\infty]=\alpha_i$. Since $\Sigma_\infty$ is $J$-holomorphic, positivity of intersections implies that all the $\alpha_i$ but one vanish, say $\alpha_1=1$, $\alpha_j=0$ $\forall j>1$. Since the classes in $\sigma_0(H_2^\hur(Y))$ have symplectic area larger than the area $a$ of the fiber by assumption, no such class can appear as a bubble, so $k=1$ and the Gromov limit consists of one component exactly.  This shows the compactness, and the rest of the proof goes as in \cite{biran1}. \vspace{,2cm}

\noindent {\it Proof of theorem \ref{t:cimo}:} We consider the situation of the theorem:
$$
(S^2,\om_\fs(a))\hra(M,\om)\overset{\pi}{\underset {\sigma_\infty}\rightleftarrows}Y,
$$
$\Sigma_\infty:=\im \sigma_\infty$ is a smooth symplectic hypersurface in $M$ whose complement $X:=M\priv \Sigma_\infty$ contains a closed Lagrangian submanifold
$L$. We fix a Riemannian metric $g$ on $L$, generic in the sense that it has only isolated closed geodesics.  \vspace{,2cm}

\begin{figure}[h!]
\begin{center} 
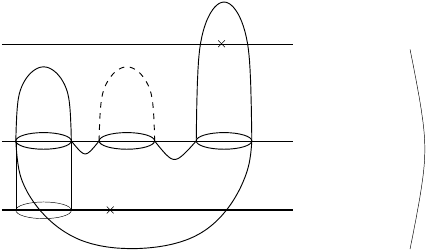
\caption{The disc $F_1\cup C_1$ in $X$ with boundary on $L$.}
\label{fig:building}
\end{center}
\end{figure}

\noindent \underline{Step 1:} (see figure \ref{fig:building}). We first work under the assumption that no closed geodesic $\gamma$ of $(L,g)$ verifies the following three conditions together:
\begin{itemize}
\its $\ell_g(\gamma)\leq \frac a\eps$ where $\eps$ is given soon,
\its $\gamma$ is contractible in $L$,
\its The Morse index of $\gamma$ is $\geq \dim L+1=n+1$. 
\end{itemize}
Let $V'\Subset V$ be two \nbds of $\Sigma_\infty$ with $L\subset M\priv V$, $p\in L$ and $\eps>0$ small enough so that there exists a symplectic embedding $$\iota:(T^*_{g,\eps}L,0_L,d\lambda_\can)\hra (M\priv V,L,\om).$$ Its image $W$ is a smoothly bounded Weinstein \nbd of $L$. Let $J$ be an $\om$-compatible almost complex structure on $V$ for which $\Sigma_\infty$ is holomorphic, and $J_n$ a sequence of  almost complex structures that coincide with $J$ in $V$ and stretch $W$'s neck along its boundary. Using lemma \ref{l:holfib}, associate to the pair $(J_n,p)$ an almost complex structure $J_n'$ that coincides with $J_n$ except on $V\priv V'$, and a pseudo-holomorphic curve 
$$
u_n:(S^2,j)\lra(M,J_n')
$$
in class $F$ and passing through $p$. By the SFT compactness theorem \cite{boelhowize}, there exists a limit holomorphic building $\cb=\{ B_i\}_{i\in \ci}$, where:
\begin{itemize}[leftmargin=0.6cm]
\its Each $B_i$ is the image of a  punctured sphere $u_i:S^2\priv \{z_1,\dots,z_{m_i}\}\to M_i$ where $M_i$ is either $W$ (the bottom layer), or $\partial W\times \R$ (the intermediate layers) or $M\priv W$ (the top layer). These layers are equipped with almost-complex structures, all denoted $\wdt J$, for which the $u_i$ are $\wdt J$-holomorphic. The important thing for us is that $\wdt J=J$ on $V$, so $\Sigma_\infty$ is a $\wdt J$-holomorphic hypersurface in the top layer, and that the holomorphicity of the $u_i$ guarantees that their symplectic areas are non-negative. 
\its $\sum [u_i]=F$, in particular each $B_i$ has area at most $a$. 
\its The map $u_i$ is asymptotic at its puncture $z_j$ to a smooth curve $\wdt \gamma_i^j\subset \partial W$ whose preimage by $\iota$ in $\partial T_{\eps,g}^\eps L$ is the lift of a geodesic $\gamma_i^j\subset L$. The orientation on $S^2$ naturally endows $\wdt \gamma_i^j$ with an orientation, that coincides with the natural orientation of the lift of $\gamma_i^j$ if $B_i\subset W$ ($z_j$ is a positive puncture) or with its opposite when $B_i\subset M\priv W$ (negative punctures). The intermediate layers have both positive and negative punctures. 
\its The different components $B_i$ of $\cb$ can be glued along the asymptotic orbits $\wdt \gamma_i^j\subset \partial W$, to form a symplectic sphere. The same holds for the components of any subbuilding $F\subset \cb$, except that the gluing may give rise to several connected components, whose domains are punctured spheres.
\its $p\in \cup B_i$. Up to reindexing, we assume that $p\in B_0$, in particular $B_0$  lies in the bottom layer $W$. 
\end{itemize}
Assume for a moment that $B_0$ has exactly one puncture. Its asymptotic $\wdt \gamma$ would be the lift of a geodesic $\gamma$ that verifies:
\begin{itemize}[leftmargin=0.6cm]
\its $\gamma$ is contractible in $L$ because $\wdt \gamma$ would be so in $T^*L$ as the boundary of the topological disc $B_0\subset W$.
\its $\ell_g(\gamma)\leq \frac a\eps$, because 
$$
a\geq \ca_\om(B_0)=\int_{\wdt \gamma} \lambda_\can=\eps\ell_g(\gamma).
$$
\its The Morse index $\mu_\gamma$ of $\gamma$  is at least $n+1=\dim L+1$. Indeed, since $p\in B_0$ and $J$ can be chosen generic, the virtual dimension of the moduli space of $J$-holomorphic planes asymptotic to $\wdt \gamma$ must be at least $2n-2$. Thus 
$$
n-3+\mu_\gamma\geq 2n-2.
$$
\end{itemize}
Since these three conditions are never met in this first step, $B_0$ has at least two punctures. 
 Let $\wdt \gamma_i$, $i=1,\dots,m_0\geq 2$ be the asymptotics of $u_0$ at its different punctures, and $\gamma_i$ the underlying geodesic loops on $L$. Since the gluing of all the components of $\cb$ is a sphere, the subbuilding $\cb\priv \{B_0\}$ has exactly $m_0$ connected components,   $F_1,\dots, F_{m_0}$, which are spheres with exactly one negative puncture, glued to $B_0$ along $\wdt \gamma_1,\dots,\wdt \gamma_{m_0}$, respectively. Since the $F_i$ are $J$-holomorphic  in $V$, they intersect the $J$-holomorphic hypersurface $\Sigma_\infty$ positively,  so exactly one of the components intersect $\Sigma_\infty$, transversally. Since these components are at least $2$, one of them, say $F_1$ lies in $X=M\priv \Sigma_\infty$. It is a standard fact that gluing to $F_1$ the cylinder $C_1$ parametrized by \fonction{\phi_1}{S^1\times [0,1]}{W}{(t,s)}{\iota(s\wdt \gamma_1(t))} provides a symplectic disc $D_1$ with boundary on $L$, therefore lying in $X$. Now 
$$
\begin{array}{ll}
\ca_\om(D_1)&=\ca_\om(F_1)+\lambda_\can(\wdt \gamma_1)\\
& \leq \ca_\om(F_1)+\lambda_\can(\wdt \gamma_1)+\dots +\lambda_\can(\wdt \gamma_{m_0})\\
& =\ca_\om(F_1)+\ca_\om(B_0)\\
&\leq  \ca_\om(\cb)=a. 
\end{array}
$$ 
In the previous chain of inequalities, the second follows from the fact that the $\lambda_\can$-actions of the loops $\wdt \gamma_i$ are positive (up to the factor $\eps$, they are the lengths of their underlying geodesics), and the third from the exactness of the symplectic form in $W$. \vspace{,2cm}

\noindent\underline{Step 2:} Now we show that the general case can be reduced to the particular setting considered in step 1. Let $A$ denote the set of geodesics of $(L,g)$ that satisfy the three forbidden conditions above. Since $g$ has isolated closed geoedesics, $A$ is a finite set. Let $\mu$ be the maximal Morse index of the geodesics in $A$.  Consider now the Lagrangian embedding 
$$
L':=L\times \T^\mu\hra (M',\om'):=(M\times \T^{2\mu},\om\oplus\om_\st). 
$$
Since the torus is aspherical, $(M',\Sigma':=\Sigma\times \T^{2\mu},L')$ verify the same properties as $(M,\Sigma,L)$, with the same constant $a$. Moreover, if we endow $L'$ with the product metric $g':=g\oplus \delta g_{\eucl}$ with $\delta$ small enough, $T_{g',\eps}^* L'$ still embeds into $M'\priv \Sigma'$. Now since the metric on the torus factor is Euclidean, a contractible geodesic of $(L',g')$ is of the form $(\gamma\times\{\star\})$, where $\gamma$ is a geodesic of $L$ and $\star\in \T^\mu$. As a result, the Morse indices and lengths of these contractible geodesics are the same as for their projections on $L$. Thus, the contractible geodesics of length  at most $\frac a\eps$ all have Morse index $\mu<\dim L'+1$. Step 1 theorefore guarantees the existence of a symplectic disc $u':(D',\partial D')\to (M',L')$ with area at most $a$. Since $\pi_2(\T^{2\mu},\T^\mu)=0$, the projection of $u'$ to the torus factor vanishes in homology, so its symplectic area vanishes. Since the symplectic form is split, the projection $u$ of $u'$ to the first factor therefore gives a disc with boundary on $L$ of area the same as $u'$, thus positive and $\leq a$. Although it is not clear whether $u:(D,\partial D)\to (X,L)$ is indeed a symplectic disc, we can argue as follows to conclude. First in dimension at least $6$, relative $h$-principle for the iso-symplectic relation in codimension $4$ implies that $u$ is homotopic to a symplectic disc, relative to its boundary, and this disc has the same area $\leq a$. In dimension $4$, the relative  $h$-principle fails, but the positivity of the area of $u$ implies that $L$ is not the sphere, so it admits a metric with non-positive curvature, which satisfies the conditions of step 1.\cqfd


Theorem \ref{t:rigskel} directly follows from the previous lemma, {\it via} a construction due to Biran \cite{biran1}. A reader unfamiliar with symplectic disc bundles may wish to consult section \ref{s:sdb} before reading the following proof. \vspace{,2cm}

\noindent{\it Proof of theorem \ref{t:rigskel}:} Let $(M^{2n},\om)$ be an integral symplectic manifold, $(\Sigma_k,\lambda_k)$ a polarization 
of $(M,\om)$ of degree $k$ with skeleton $\Gamma_k$ and let $L\subset M\priv \Gamma_k$ a closed Lagrangian submanifold. Biran decomposition theorem shows that $L$ embeds into the symplectic disc bundle $\sdb(\Sigma_k,\nf 1k)$ of width $\nf 1k$. Since $L$ is compact, it even lies in a sub-disc bundle of width $\frac 1k-\eps$, that compactifies into a ruled symplectic manifold $(Y,\om)$ of width $\nf 1k-\eps$ that contains two disjoint symplectic sections $\sigma_{0,\infty}:\Sigma_k\hra Y$.  Here, $\im \sigma_\infty$ corresponds to the compactification of the boundary of our symplectic disc bundle, so  has image disjoint from $L$, and $\sigma_0^*\om=\om_{|\Sigma_k}$, thus has integral class. The conditions in  theorem \ref{t:cimo} are therefore met, so we get a symplectic disc $(D,\partial D)\subset (Y,L)$ of area $<\nf 1k-\eps$ that does not intersect $\Sigma_\infty$. This disc therefore belongs to $\sdb(\Sigma_k,\nf 1k)\subset M\priv \Gamma_k$, and fulfills our requirements.  \cqfd

Notice that if in theorem \ref{t:rigskel}, we also assume that $L$ lies in the complement of a symplectic divisor with normal crossings $\ch\subset M$ which is invariant under the Liouville flow of $\lambda_k$, then $L$ lies in $\pi^{-1}(\Sigma_k\priv \ch)\subset \sdb(\Sigma_k,\nf 1k)$ (where $\pi:\sdb(\Sigma_k,\nf 1k)\to \Sigma_k$ is the natural projection). In the compactification $(Y,\om)$ defined in the previous proof, $\pi^{-1}(\ch)\cup \Sigma_\infty$ is a symplectic divisor with normal crossing, which can be made $J$-holomorphic for some $\om$-compatible almost-complex structure $J$. Then the proof of theorem \ref{t:cimo} can be carried on, and using positivity of intersection between holomorphic curves and $J$-holomorphic divisors as before, it provides a symplectic disc
$$
u:(D,\partial D)\to \big(Y\priv(\Sigma_\infty\cup \pi^{-1}(\ch)),L\big)\approx \big(M\priv (\ch\cup \Gamma_k),L\big). 
$$
We therefore get the following relative version of theorem \ref{t:rigskel}:
\begin{theorem}\label{t:rigskelaff} 
Let $(M,\om)$ be a symplectic manifold with integral symplectic class, $\ch\subset M$ a symplectic divisor with normal crossings and $X:=M\priv \ch$. Let also  $(\Sigma_k,\lambda_k)$ a smooth polarization of $(M,\om)$ of degree $k$ whose Liouville flow preserves $\ch$, and $\Gamma_k$ its skeleton. Then any closed Lagrangian submanifold $L$ in $X\priv \Gamma_k$ bounds a symplectic disc in $X\priv \Gamma_k$ of area $<\nf 1k$. In particular:
$$
c_\lag(X\priv \Gamma_k)\leq \frac 1k. 
$$ 
\end{theorem}
If $(M,\om)=(\P^n,\om_\fs)$ or $(\CP^1,\om_\fs)^n$ and $\ch=\{z_0=0\}$ or $\ch=\cup\{z_i=+\infty\}$, respectively, it is rather easy to produce polarizations of any degree whose Liouville flow preserves $\ch$. The last result therefore estimates the Lagrangian capacity of the complement of the skeleton in the ball or the polydisc, respectively (see the end of section \ref{s:examples}).

\subsection{Lagrangian lift of the sweep of a Legendrian submanifold under a Hamiltonian flow.}
The following proposition is a reformulation of Mohnke's trick \cite{mohnke}, that we already used in a slightly weaker form in \cite{opsc}.
\begin{proposition}\label{p:laglift}
Let $(M,\alpha)$ be a contact manifold, $\Lambda\subset M$ be a closed Legendrian submanifold and $H\in \cc^\infty(M\times[0,T])$ a smooth positive time-dependent contact Hamiltonian function  such that the map \fonction{\Psi}{\Lambda\times[0,T]}{M}{(x,s)}{\Phi^s_H(x)} is an embedding. 
Then there exists a Lagrangian embedding $L$ of $\Lambda\times S^1$ in the symplectization $SM:=(M\times \R_+^*,d(R\alpha))$  with the following properties:
\begin{itemize}
\its if $\pi:SM\to M$ denotes the obvious projection, $\pi(L)$ is the sweep of $\Lambda$ under the flow of $H$:
$$
\pi(L)=\{\Phi_H^s(x),\; x\in \Lambda, s\in [0,T]\}.
$$
In particular, if $H$ has no chord between $\Lambda$ and some closed subset $K$, then $L\cap (K\times\R_+^*)=\emptyset$. 
\its $\ca_{\min}(L)=T$.
\its If $H\geq 1$ then $L\subset S^-M:=M\times (0,1]$. 
\end{itemize} 
\end{proposition}
\noindent {\it Proof:} We recall that the contact isotopy associated to $(H,\alpha)$ is the flow of the time-dependent vector field $X_s$ defined by 
$$
\left\{
\begin{array}{l}
\alpha(X_s)=H_s,\\
d_x\alpha(X_s,\cdot)=R_\alpha\cdot H_s(x)\alpha-dH_s.
\end{array}\right.
$$
 Since $\Lambda$ and $K$ are closed and $H$ is smooth, $H$ can be extended to $M\times [0,T+\eps]$ for some small $\eps>0$, still with no Reeb-chord from $\Lambda$ to itself or to $K$.   
Define  the map \fonction{\Phi}{\Lambda\times [0,T+\eps]\times (0,1]}{M\times (0,+\infty)}{(x,s,t)}{(\Phi_H^s(x),\nf t{H_s(x)}),}
which is an embedding because of the injectivity of $\Psi$. 
Then 
$$
\begin{array}{lll}
&\Phi^*\alpha\big(\frac \partial{\partial s}\big)&=\alpha(X_s)=H_s(x),\\
&\Phi^*\alpha\big(\frac \partial{\partial t}\big) & =0,\\
\forall u\in T\Lambda,& \Phi^*\alpha(u)& =\alpha(\Phi_*u)=\alpha(\Phi_{H*}^su)=0.
\end{array}
$$
The vanishing of the last expression holds because $\Phi_H^s$ is a contact diffeomorphism and $u\in T\Lambda\subset \ker \alpha$, so $\Phi_{H*}^su$ remains in the contact distribution $\ker \alpha$ for all $s$. 
Thus
$$
\Phi^*\alpha=H_sds 
$$
and
$$
\Phi^*(R\alpha)  =\frac t{H_s}\Phi^*\alpha =tds.
$$
Let now $\gamma=(s(\theta),t(\theta)):S^1\to [0,T+\eps]\times (0,1]$ be an embedded loop that encloses an area $T$ (for the form $dt\wedge ds$) and 
$$
L:=\Phi(\Lambda\times \im\gamma)=\{(\Phi_H^{s(\theta)}(x),\frac {t(\theta)}{H_{s(\theta)}(x)}),x\in \Lambda,\theta\in S^1\}. 
$$
By the above computation, $L$ is an embedded Lagrangian submanifold, that obviously projects to the $H$-sweep of $\Lambda$. If $H\geq 1$, $L\subset\{R\leq 1\}=S^-M$. Finally, since the symplectic form on $SM$ is $d(R\alpha)$, $\ca_{\min}(L)$ is simply the minimal positive $R\alpha$-action of the loops in $L$. Since $H_1(L)\overset {\Phi_*}\simeq H_1(\Lambda)\oplus \langle [\gamma]\rangle$ and $\Phi^*(R\alpha)=tds$ we see that the area morphism vanishes on $H_1(\Lambda)$ and equals $T$ on $[\gamma]$. Thus indeed, 
$\ca_{\min}(L)=T$.\cqfd \vspace{-,3cm}

\subsection{On the area class of  skeleta.} 
Recall the definition of polarization:
\begin{definition}\label{d:weincomp} 
Let $(M,\om)$ be a closed symplectic manifold with integral cohomology class. A polarization of degree $k$ is a 
triple $(\Sigma,\lambda,\phi)$ composed of a symplectic hypersurface  Poincaré-dual to $k[\om]$, a Liouville form $\lambda$ on $M\priv \Sigma$ and a Morse function $\phi$ on $M\priv \Sigma$ that tends to $-\infty$ at $\Sigma$, such that the Liouville vector field $X_\lambda$ is gradient-like for $\phi$. 

The skeleton of $(\Sigma,\lambda)$ is defined as the set of points that do not reach $\Sigma$ under the Liouville flow of $\lambda$ and is denoted $\skel(\Sigma,\lambda)$. It  is an isotropic CW-complex, composed of the stable manifolds of the critical points of $\phi$. 
\end{definition}
 Since the skeleton of a polarization is isotropic, we can define an area morphism \fonction{\ca}{H_2(M,\skel(\Sigma,\lambda),\R))}{\R}{A}{[\om]\cdot A}
\begin{proposition}\label{p:exskel}
Let $(M,\om)$ be an integral symplectic manifold with a polarization $(\Sigma,\lambda)$ of degree $k$. Then the image of the area morphism of $\skel(\Sigma,\lambda)$ lies in $\frac 1k\Z$. 
\end{proposition}
\noindent {\it Proof:} The Liouville form $\lambda$ has residue $-\frac 1k$ along $\Sigma$, in the sense that its integral along a loop that turns once and positively around $\Sigma$ is very close to $-\frac 1k$ when this loop is $\cc^1$-close to a constant in $\Sigma$. This can be shown by considering any $2$-cycle $S$ in $\Sigma$, perturbing it to a $2$-cycle $\wdt S$ in $M$ that intersects $\Sigma$ transversally, and computing its area in two different ways, one using the Poincaré-Duality condition of $\Sigma$ and  one   using Stokes theorem on the complement of very small discs around the intersections of $\wdt S$ with $\Sigma$ (see for instance \cite{opshtein3} for a complete proof). 
 
Let now $S$ denote a $2$-cycle in $M$ with boundary on $\skel(\Sigma,\lambda)$, and $\partial S$ the boundary cycle in $\skel(\Sigma,\lambda)$. 
Since $\skel(\Sigma,\lambda)$ is the disjoint union of the stable manifolds of the critical points of $\phi$, we can write 
$$
\partial S=\cup_{x\in \crit \phi} \gamma_x,
$$
with $\gamma_x:=\partial S\cap W^s(x)$. Under the positive Liouville flow (that leaves $\skel(\Sigma,\lambda)$ invariant), $\gamma_x$ deforms to a union of Liouville trajectories that connect $x$ to other critical points in the closure of $W^s(x)$. We therefore get a new representative $\gamma$ of $[\partial S]\in H_1(\skel(\Sigma,\lambda))$ made of a union of orbits of the Liouville flow. Since $\skel(\Sigma,\lambda)$ is isotropic and $\lambda$ is a Liouville form, the $\lambda$-action of $\partial S$ coincides with that of $\gamma$. Now since $\gamma$ is tangent to $X_\lambda$, on which $\lambda$ vanishes (because $\lambda(X_\lambda)=\om(X_\lambda,X_\lambda)=0$), 
$$ 
\int_\gamma \lambda=0. 
$$
In view of the preliminary remark on the residue of $\lambda$, we therefore see that 
$$
\ca_\om(S)=\frac 1k[S]\cdot [\Sigma] \in \frac 1k\Z
$$
(in this equality, the intersection means homological intersection between a $2$-cycle relative to $\skel(\Sigma,\lambda)$ and a $2n-2$-cycle that is disjoint from $\skel(M,\lambda)$). \cqfd

As a consequence, we get the first assertion of theorem \ref{t:legbarp}. 
\begin{corollary}\label{c:legbardot1}
Let $(N^{2n},\tau)$ be a closed symplectic manifold with integral cohomology class and $(\cp^{2n+1},\alpha)\overset \pi\to (N,\tau)$ its prequantization contact bundle. Let $(\Sigma,\lambda)$ be a polarization of $(N,\tau)$ of degree $k$, and $\Gamma_k$ its skeleton. Let $x\in \Gamma_k$ and $e:=\{e_1,\dots,e_k\}\in \pi^{-1}(x)$ be $k$ equi-distributed points in the fiber of $x$ (with respect to the $S^1=\R_{/\Z}$-parametrization).  There exists a unique compact Legendrian  CW-complex $\Lambda_k$ whose intersection with $\pi^{-1}(x)$ is $e$ and such that $\pi(\Lambda_k)=\Gamma_k$. We call it a $k$-fold Legendrian cover of $\Gamma_k$. 
 \end{corollary}
\noindent{\it Proof:} We recall that 
$\cp(N,\tau)\to N$ is a circle bundle on which $\alpha$ is a connection form with curvature $-\tau$. The horizontal distribution is $\ker \alpha$, which is the contact distribution. Several choices for $\alpha$ are possible (for instance by adding the pull-back of a closed $1$-form on $N$), but they lead to the same contact forms, up to isotopy. Denoting $P^\gamma_\alpha$ the parallel transport along a path $\gamma:[0,1]\to N$, we have the following straightforward formula for $\alpha':=\alpha+\pi^*\vartheta$, $d\vartheta=0$:
$$
\forall e\in \pi^{-1}(\gamma(0)), \hspace{,3cm} P^\gamma_{\alpha'}(e)=e^{i\int_\gamma \vartheta}P_\alpha^\gamma(e). 
$$

Let now $G\subset H_1(\Gamma_k)$ denote the image of the boundary operator $\partial:H_2(\cp,\Gamma_k)\to H_1(\Gamma_k)$,  $F$ a complement of $G$ in $H_1(\Gamma_k)$ and $(f_1,\dots,f_n)$ a basis of $F$ (so the images of the $f_i$ in $H_1(\cp)$ give an independent family). Fix a closed $1$-form $\vartheta$ on $N$ whose values on fixed representatives $\gamma_i$ of $f_i$ 
verify $\ds P_\alpha^{\gamma_i}=e^{i\int_{\gamma_i} \vartheta}$.
Then defining $\alpha':=\alpha-\pi^*\vartheta$, we get by the formula above:
$$
P_{\alpha'}^{\gamma_i}=\id. 
$$
Since $\Gamma_k$ is isotropic, the restriction of $P_{\alpha'}^{\,\centerdot}$ to $\Gamma_k$ depends only on the homology class of the loop on which it is evaluated (in $H_1(\Gamma_k)$), so we can see $P_{\alpha'}|_{\Gamma_k}$ as a group homomorphism from $H_1(\Gamma_k)$ to $U(1)$, and we have ensured by the choice of $\alpha'$
that this morphism vanishes on our subspace $F\subset H_1(\Gamma_k)$. Moreover, since $d\alpha'=\pi^*\tau$, proposition \ref{p:exskel} guarantees that $P_{\alpha'}^{\, \centerdot}$ is a rotation by an integer multiple of $\frac 1k$ for each homology class in $G$, so we finally see that the image of $P_{\alpha'}$ in $U(1)$ is contained in the subgroup $\Z_k\subset U(1)$ of order $k$.

Let $(x,e:=\{e_1,\dots,e_k\})$ be the point in $\Gamma_k$ and the equi-distributed points in its fiber considered in the statement. 
For each $y\in \Gamma_k$, fix a path $\gamma_{xy}$ joining $x$ to $y$ in $\Gamma_k$, and define 
$$
\Lambda_k:=\{P^\alpha_{\gamma_{x,y}}(e_i), \; i\in [1,k],\; y\in \Gamma_k\}. 
$$
Since parallel transport is by isometries of $S^1$, $\Lambda_k\cap \pi^{-1}(y)$ is made of exactly $k$ equi-distributed points in $S^1=\pi^{-1}(y)$. Since $\im P_{\alpha'}\subset \Z_k$, that stabilizes any set of $k$ equi-distributed points,  this definition of $\Lambda_k$ does not depend on the choice of the system of paths $\gamma_{x,y}$.  This implies by standard arguments that $\Lambda_k$ is a compact CW-complex. It is  Legendrian since, by definition of parallel transport, the tangent spaces to $\Lambda_k$ are horizontal, so they lie in the contact distribution. \cqfd\vspace{-,3cm}

\subsection{On quasi-holomorphic sections.}\label{s:qh} In our context, quasi-holomorphic sections were introduced by Giroux in \cite{giroux1,giroux2}. We recall now their definition and relevant properties. 
Let $(V,\om,J,g)$ be a linear symplectic space with an $\om$-compatible linear complex structure $J$, that altogether define the scalar product $g=\om(\cdot,J\cdot)$. The metric $g$ provides a dual  norm on $V^*$, $\om$ an isomorphism $X_\centerdot:V^*\to V$ and  the complex structure  an isomorphism ${}^c:V^*\to V^*$ given by $\alpha^c=\alpha(J\cdot)$ (it is an anti-involution).  In fact, 
the dual norm is euclidean, $\alpha^c(X_\alpha)=|\alpha|^2$ and $X_\centerdot$ is an isometry. The last point follows from:
$$
|\alpha|=\max_{|X|=1}\alpha(X)=\max_{|X|=1}\om(X_\alpha,X)=\max_{|X|=1}g(X_\alpha,JX)=|X_\alpha|.
$$
One checks without difficulty that $X_{\alpha^c}=-JX_\alpha$. For a complex-valued linear form $\Lambda\in V^*_\C:=V^*\otimes \C$, we define 
$$
\begin{array}{l}
\Lambda^{(0,1)}:=\frac 12(\Lambda+i\Lambda^c)\\
\Lambda^{(1,0)}:=\frac 12(\Lambda-i\Lambda^c),
\end{array}
$$
so that 
$$
\Lambda^{(0,1)}(Jv)=-i\Lambda^{(0,1)}v,\hspace{,3cm} \Lambda^{(1,0)}(Jv)=i\Lambda^{(1,0)}v,\hspace{,3cm} \Lambda=\Lambda^{(0,1)}+\Lambda^{(1,0)}. 
$$
\begin{lemma}\label{l:linqh} Let $\Lambda=\alpha+i\beta\in V^*_\C$, $\alpha,\beta\in V^*$. Assume that there exists $\kappa<1$ such that 
$$
|\Lambda^{(0,1)}|\leq \kappa |\Lambda^{(1,0)}|
$$ 
(we say that $\Lambda$ is $\kappa$-quasi-holomorphic).  Then 
$$
\om(X_\beta,X_\alpha)\geq \frac{1-\kappa^2}{1+\kappa^2}\frac{|\alpha|^2+|\beta|^2}2
$$
and 
$$
|X_\alpha+JX_\beta|^2\leq \frac{2\kappa^2}{1+\kappa^2}(|X_\alpha|^2+|X_\beta|^2).
$$
In particular, the vanishing of one of either $\alpha$ or $\beta$ implies the vanishing of the other. 
\end{lemma}
\noindent{\it Proof:}
Notice first that given $\alpha\in V^*$, 
$$
|\alpha+i\alpha^c|=\sqrt 2|\alpha|: 
$$
since $J$ is an isometry of $V$, $|\alpha^c|=|\alpha|$, so  $|\alpha+i\alpha^c|\leq \sqrt 2|\alpha|$. The reverse inequality is obtained by feeding $\alpha+i\alpha^c$ with the vector $v+Jv$, where $\alpha(v)=|\alpha|$, $|v|=1$ (such a vector automatically verifies $\alpha(Jv)=0$ because $v\perp Jv$). Now 
$$
\begin{array}{l}
2\Lambda^{(1,0)}=\Lambda-i\Lambda^c=\alpha+\beta^c-i(\alpha+\beta^c)^c\\
2\Lambda^{(0,1)}=\Lambda+i\Lambda^c=\alpha-\beta^c+i(\alpha-\beta^c)^c\\
\end{array}
$$
By the previous remark, the $\kappa$-quasi-holomorphicity of $\Lambda$ gives 
$$
|\alpha-\beta^c|\leq \kappa|\alpha+\beta^c|. 
$$
Squaring the last inequality and developping gives 
$$
g(\alpha,\beta^c)\geq \frac{1-\kappa^2}{1+\kappa^2}\frac{|\alpha|^2+|\beta|^2}{2}
$$
so
$$
\om(X_\beta,X_\alpha)=
g(X_\alpha,-JX_\beta)=g(X_\alpha,X_{\beta^c})=g(\alpha,\beta^c)\geq \frac{1-\kappa^2}{1+\kappa^2}\frac{|\alpha|^2+|\beta|^2}{2}.
$$
Then,
$$
|X_\alpha+JX_\beta|^2=|X_\alpha|^2+|X_\beta|^2-2g(X_\alpha,-JX_\beta)\leq \frac{2\kappa^2}{1+\kappa^2}(|X_\alpha|^2+|X_\beta|^2).\hspace{,2cm}\square
$$
\vspace{,1cm}

Given now a symplectic manifold endowed with an almost complex structure $(M,\om,J,g)$ and a Hermitian line bundle $\cl\to M$, the definitions above obviously apply pointwise, and provide a map $X_\centerdot:\Om^1(M,\R)=\Gamma(T^*M)\to \Gamma(TM)$, dual norms on $TM^*\otimes \cl$, and fiberwise isomorphisms $^c:T^*M\otimes \cl\to T^*M\otimes \cl$. Given a Hermitian connection $\nabla$ on $\cl$, we can define its linear and anti-linear parts:
$$
\nabla':=\nabla^{(1,0)} \hspace{,5cm}
\nabla'':=\nabla^{(0,1)}. 
$$
A section $s:M\to \cl$ is said  $\kappa$-quasi-holomorphic  if 
$$
\forall x\in M, \hspace{,3cm}|\nabla'' s(x)|\leq \kappa |\nabla's(x)|. 
$$
\begin{lemma}\label{l:qhsec} Let $(M,\om,J,g)$ be a symplectic manifold with a compatible almost complex structure, $\cl\to M$ a Hermitian line bundle over $M$ and $s:M\to  \cl$ a $\kappa$-quasi-holomorphic section, with $\kappa<1$. On $\{s\neq 0\}$, $\nabla s=(d\ln|s|-i\beta)s$, $\beta\in \Om^1(M,\R)$, and 
\begin{itemize}
\its 
$X_\beta\cdot (-\ln |s|)\geq \frac {1-\kappa^2}{1+\kappa^2}\frac{|\alpha|^2+|\beta|^2}{2}$,
\its  $\nabla s(p)=0 \Longleftrightarrow d\ln|s|(p)=0$,
\its  $\nabla s(p)$ either vanishes or is surjective.  
\end{itemize}
\end{lemma}
\noindent{\it Proof:} 
On $\{s\neq 0\}$, we can write for any smooth complex-valued function $f$: 
$$
\nabla f\frac s{|s|}=(df +f\vartheta)\frac s{|s|},
$$
and $\vartheta$ is a connection $1$-form. Since the connection is Hermitian, $\vartheta=-i\beta$, $\beta\in \Om^1(M,\R)$. Applying the formula for $f=|s|$ further gives
$$
\nabla s=(d|s|-i|s|\beta)\frac  s{|s|}=(d\ln|s|-i\beta)s. 
$$
This concludes the first bullet. 
The remaining assertions just follow by applying lemma \ref{l:linqh} pointwise. \cqfd

\subsection{Connections on powers of line bundles}
Along this paper, sections of complex line bundles of the form $\cl^{\otimes k}$ (denoted for short $\cl^k$)  will be used intensively. These power line bundles are naturally equipped with connections as long as $\cl$ is, and the connections we consider throughout the whole paper are of this kind. This is very standard, but we recall it for self-containedness, highlighting the points of further use. 

\begin{proposition}\label{p:nablak} 
Let $(\cl\overset \pi\lra M,\nabla)$ be a hermitian line bundle equipped with a unitary connection. There is a natural unitary connection $\nabla^k$ on $\cl^{ k}$ defined in the following way:
\begin{itemize}
\its On an open set $U\subset M$ above which $\cl$ is trivial, pick a non-vanishing section $e:U\to \cl$ and write $\nabla e=i\theta e$, $\theta\in \Om^1(M)$ (the form is real-valued because $\nabla$ is unitary),
\its Define $\nabla^ke^{\otimes k}:=k\theta e^{\otimes k}$,
\its Extend $\nabla^k$ on all sections above $U$ by the Leibniz rule.    
\end{itemize}
This connection has the following properties:
\begin{enumerate}
\item If $s:[0,1]\to \cl$ is a $\nabla$-parallel section along $\gamma:[0,1]\to M$ then $s^{\otimes k}$ is $\nabla^k$-parallel along $\gamma$. 
\item If $M$ is equipped with an almost complex structure $J$ and $s:M\to \cl$ is $\kappa$-quasi-holomorphic, so is $s^{\otimes k}$. 
\item If $\nabla$ has curvature $\vartheta$, $\nabla^{\otimes k}$ has curvature $k\vartheta$. 
\end{enumerate}
\end{proposition}
\noindent {\it Proof:} The three consequences are completely straightforward and we leave them to the reader. Checking that $\nabla^k$ indeed provides a well-defined connection on $\cl^k$ is equivalent to the fact that it is natural in the sense that it does not depend on the section $e$ we pick in this definition. We prove this point. Let $(e,\theta,\nabla^k)$ be defined on $U$ as proposed in the statement. Let $e'$ be another non-vanishing section above $U$ and $\theta'\in \Om^1(U)$ its associated $1$-form, meaning that $\nabla e'=\theta' e'$. Then $e'=fe$  for some smooth function $f:U\to \C^*$ and $\nabla e'=(\frac{df}f+\theta)e'$, so $\theta'=\theta+\frac{df}f$. On the other hand, $e'{}^{\otimes k}=f^ke^{\otimes k}$ so 
$$
\nabla^k e'{}^{\otimes k}=(df^k+k\theta)e^{\otimes k}=k(\frac{df}f+\theta)e'^{\otimes k}=k\theta'e'{}^{\otimes k}. 
$$  
 This shows indeed that this definition of $\nabla^k$ does not depend on the chosen section $e$. \cqfd

\vspace{-,3cm}
\begin{remark}
The three properties listed above have implications of different nature. The point ($3$) tells that these connections are compatible with the symplectic category, in the sense that if $\nabla$ has curvature $-i\om$, $\nabla^k$ has curvature $-ik\om$. The point ($2$) shows that these connections are compatible with the complex category, allowing to compute examples in the K\"ahler setting, for instance the holomorphic $\nabla^k$ sections on $\co(k)$ are holomorphic in the usual sense provided the connection $\nabla$ on $\co(1)$ is the most natural connection. Finally ($1$) happens to be technically very usefull in our proof. 
\end{remark}

\section{Negative symplectization of prequantization bundles.}\label{s:sdb}
In this section we introduce symplectic disc bundles and discuss their relations with prequantization bundles. Most of it is very classical and already explained in Biran's orginal paper, where symplectic disc bundles first appear \cite{biran1}. Sections \ref{s:blowup} and \ref{s:spb} are a bit more original, describing a rather obvious but non-symplectic compactification of these objects. Leaving the purely symplectic world has a price in terms of the attention we need to pay, but provides the sharpest possible results. Degeneracy is dealt with an unusual almost-complex structure, introduced in section \ref{s:sdbdef}. 

\subsection{Symplectic disc bundles.}\label{s:sdbdef}
Let $(N,\tau)$ be a  symplectic manifold with integral symplectic class and  $\ell\in \N^*$. Let $\cl\overset \pi\to (N,\tau)$ be a line bundle with $c_1(\cl)=\ell[\tau]$
equipped with a radial coordinate $R$ and a connexion $1$-form $\Theta$ that satisfies
$$
\left\{
\begin{array}{l}
\Theta|_{\cl_p}=d\theta,\\
d\Theta=-\ell \pi^*\tau.
\end{array}\right.
$$
(Here $p\in N$ and $\cl_p$ denotes its fiber in $\cl$). 
The closed two-form 
\begin{equation}\label{e:omsdb}
\om_0:=\pi^*\tau+d(R\Theta)=(1-\ell R)\pi^*\tau+dR\wedge \Theta
\end{equation}
is symplectic on $\{R<\frac 1\ell\}$.  This subset of $\cl$ is the symplectic disc bundle over $(N,\tau)$ associated to $\cl$:
$$
\sdb(N,\tau,\ell):=\big(\{R<\frac 1\ell\}\subset \cl,\om_0\big). 
$$
The $1$-form
$$
\lambda_0:=(R-\frac 1\ell)\Theta,
$$
defined on the complement of the zero section $N_0:=0_\cl$, is a Liouville form for $\om_0$ with radial Liouville vector field
$$
X_0:=(R-\frac 1\ell)\frac\partial{\partial R}. 
$$
One main feature of these symplectic disc bundles is that they provide explicit local models for some symplectic divisors:
\begin{proposition}\label{p:locsdb} 
Let $(M,\om)$ be a symplectic manifold and $\Sigma\subset M$ a smooth symplectic divisor. If $[\Sigma]=\pd(\ell [\om])$, or more generally if  the symplectic normal bundle to $\Sigma$ has first Chern class $\ell[\om_{|\Sigma}]$, then $\op(\Sigma,M)$ can be symplectically identified with 
$$
\sdb_\eps(\Sigma,\om_{|\Sigma},\ell):=\{R<\eps\}\subset \sdb(\Sigma,\om_{|\Sigma},\ell)
$$
for some $\eps\ll 1$. 
\end{proposition}

In most of the present paper, we only need to consider 
$$
\sdb(N,\tau):=\sdb(N,\tau,\ell=1). 
$$ 
 Given a $\tau$-compatible almost complex structure $j$ on $N$, one may lift $j$ to $\sdb(N,\tau)$ by pulling back the natural almost complex structure on $\cl$ by the  radial map 
$$
(z,R,\theta)\longmapsto (z,\frac R{1- R},\theta),
$$
that takes the disc bundle to the whole of $\cl$. On checks that 
$$
\left\{
\begin{array}{l}
 J(T\cl_p)= T\cl_p \text{ and } J(\ker \Theta)=\ker \Theta,\\
 J \frac\partial{\partial R}=\frac 1{2R(1-R)}\frac\partial{\partial \theta},\\
 J(U)=\pi^*j\pi_*U\hspace*{,3cm} \forall U\in \ker \Theta
\end{array}\right.
$$
(in the last line and throughout the paper, $\pi^*$ denotes the inverse of the restricted isomorphism $\pi_{*|\ker \Theta}:\ker \Theta\to TN$ when it acts on a vector of $TN$).

\subsection{Blow-up of the zero-section, contact compactification.}\label{s:blowup}
The negative part of the symplectization of the prequantization bundle $(\cp,\alpha)\to (N,\tau)$, denoted $S^{\leq 0}\cp$ is defined by 
$$
S^{\leq 0}\cp:=(\cp\times (0,1],d(R\alpha)). 
$$
It is an open symplectic manifold, with two ends, for $R\in \{0,1\}$. The end  $\{R=1\}$ is a contact type boundary, with periodic Reeb vector field, so it  can  be compactified by standard reduction procedure. This compactification leads to the symplectic disc bundle over $(N,\tau)$: 
\begin{lemma}\label{l:sdbsp}
Let $(N,\tau)$ be a symplectic manifold with integral symplectic class, $(\cp,\alpha)$ its prequantization bundle and $\sdb(N,\tau)$ its symplectic disc bundle. The map \fonction{\Phi}{\big(\sdb(N,\tau)\priv N_0,\lambda_0\big)}{S^{<0}\cp=(\cp\times (0,1),R\alpha)}{(z,R,\theta)}{(z,-\theta,1-R)} is a well-defined exact symplectomorphism ($\Phi^*(R\alpha)=\lambda_0$).
\end{lemma}
\noindent {\it Proof:} The map $\Phi$ is well-defined because $S\cp$ is a complex line bundle over $N$ of Chern-class equals to the Euler class of $\cp\to N$, which is $[\tau]$, hence opposite to the Chern-class of $\cl$. Obviously, 
$\Phi^*(R\alpha)=(1-R)\Phi^*\alpha$ and since $\alpha$ is a connection $1$-form on $\cp$ with curvature $\tau$, $\Phi^*\alpha$ is a $1$-form on $\sdb(N,\tau)$, invariant under radial transformations, and 
$$
\begin{array}{l}
\Phi^*\alpha|_{F}=-d\theta\\
\Phi^*d\alpha=\pi^*\tau.
\end{array}
$$
Thus $-\Phi^*\alpha$ is a connection form on $\sdb(N,\tau)\priv N_0$ with curvature $-\pi^*\tau$, so we can assume that it coincides with $\Theta$. Thus, 
$$
\Phi^*(R\alpha)=(1-R)\Phi^*\alpha_0=(R-1)\Theta=\lambda_0. \hspace*{2cm}\square
$$

\begin{lemma}\label{l:skelift}
Let $(\cp,\alpha)\to (N,\tau)$ be a contact prequantization bundle. Let $(\Sigma,\lambda)$ be a polarization of $N$ of degree $k$ and $\Gamma$ its skeleton (see definition \ref{d:weincomp}). Let $\Lambda$ be the $k$-fold Legendrian lift of $\Gamma$ to $\cp$ described in corollary \ref{c:legbardot1} and $\Phi:\sdb(N,\tau)\priv N_0\to S^{<0}\cp$ the symplectic identification provided by lemma \ref{l:sdbsp}. Then 
$$
\Phi^{-1}(\Lambda\times (0,1))
$$
is the unique $\om_0$-isotropic CW-complex that  projects under $\pi$ to $\Gamma$, and such that the fibers of this projection are $k$ radial rays. This uniqueness is meant up to global fibered rotation of fixed angles, {\it i.e.} maps of the form $\cR_c:(z,R,\theta)\mapsto (z,R,\theta+c)$.  
\end{lemma}
\noindent{\it Proof:} The subset $\Phi^{-1}(\Lambda\times (0,1))$ is indeed an isotropic CW-complex and its projection to $N$ is 
$$
\pi\circ\Phi^{-1}(\Lambda\times(0,1))=\pi(\Lambda\times (0,1))=\pi(\Lambda)=\Gamma. 
$$
Since $\pi_{|\Lambda}$ is a $k$-fold cover of $\Lambda$ and since $\Phi^{-1}(\{\star\}\times(0,1))$ is a radial ray, the fibers of $\pi$ in $\Phi^{-1}(\Lambda\times(0,1))$ are indeed $k$ radial rays. 
Further, any global rotation of this subset (by a map $\cR_c:(z,R,\theta)\mapsto (z,R,\theta+c)$ where $c$ is a constant) is of the form 
$\Phi^{-1}\big(\Lambda'\times (0,1)\big)$, where $\Lambda'$ is obtained from $\Lambda$ by a fibered rotation in $\cp$ with constant angle. Thus $\Lambda'$ is a Legendrian $k$-fold cover of $\Gamma$, so $\pi(\Phi^{-1}(\Lambda'\times (0,1)))=\Gamma$ as well. It remains to check the uniqueness part of the statement. 

Let $L$ be an isotropic radially invariant CW-complex of $\sdb(N,\tau)\priv N_0$ whose projection to $N$ is a $k$-fold cover of $\Gamma$ (in the sense of the lemma). Then $\Phi(L)$ is isotropic and radially invariant in $S^{<0}\cp$, so is of the form $\Lambda'\times(0,1)$ for some Legendrian CW-complex $\Lambda'$.   Since moreover 
$$
\Gamma=\pi(L)=
\pi(\Lambda'\times (0,1))=\pi(\Lambda'),
$$
$\Lambda'$ is  a Legendrian $k$-fold lift of $\Gamma$. By corollary \ref{c:legbardot1}, $\Lambda'$ is obtained from $\Lambda$ by a global rotation, so $\Phi^{-1}(\Lambda'\times (0,1))$ is obtained from $\Phi^{-1}(\Lambda\times (0,1))$ by a global rotation as well. \cqfd
\vspace*{-,3cm}
\subsection{Projectivization of symplectic disc bundles.}\label{s:spb} The map $\Phi$ provided by lemma \ref{l:sdbsp} allows to compactify  the negative symplectization  $S^{<0}\cp(N,\tau)$ of $\cp(N,\tau)$  at its end  $\{R=1\} $ by gluing $N_0$. 
We now aim at compactifying the end $\{R=0\}$ of  $S^{<0}\cp$, which can be obtained equivalently as a compactification of $\sdb(N,\tau)$ at the boundary $\{R=1\}$ of the disc bundle, still  by lemma \ref{l:sdbsp}. In \cite{biran2}, Biran used a compactification of $\sdb(N,\tau)$ obtained in the following way:
\begin{itemize}[leftmargin=*]
\its Remove a \nbd  of the boundary, thus considering $\{R<{1-\delta}\}\subset \sdb(N,\tau)$,
\its Consider $\delta\sdb(N,\tau,-1):=\big(\sdb(N,\tau,-1),\delta\om_0\big)$.
\its Observe that $\{0<R<\delta\}\subset \delta\sdb(N,\tau,-1)$ is exact symplectomorphic to $\{1-2\delta<R<1-\delta\}\subset \sdb(N,\tau)$, {\it via} the most simple map  $\psi:(z,R,\theta)\mapsto (z,1-\delta-R,-\theta)$. 
\its Define 
$$
\spb(N,\tau,\delta):=\{R<1-\delta\}_{\subset\, {\sdb(N,\tau)}}\underset \psi\bigcup \{R<\delta\}_{\subset \,\delta\sdb(N,\tau,-1))}.
$$ 
\end{itemize}
The result of this compactification is a symplectic sphere bundle over $N$ (hence a ruled symplectic manifold), the $0$ and  $\infty$-sections being symplectic, namely $(N,\tau)$ and $(N,\delta\tau)$, respectively. The area of the fibers is $1-\delta$, and this compactification contains $\sdb_{1-\delta}(N,\tau)$. 
Unfortunately, this symplectic manifold has symplectic class in $\frac 1{N(\delta)}\Z$, with $N(\delta)\to +\infty$ when $\delta\to 0$, so has polarizations only in very large degrees when $\delta$ becomes small, which makes it useless in this paper. 

Instead, we describe now a non-symplectic compactification of $\sdb(N,\tau)$ into a sphere bundle with a smooth closed $2$-form that degenerates along the section at infinity.  As we have already seen, $\sdb(N,\tau)$ lies in a hermitian line bundle $\cl$, on which we have a smooth $2$-form 
$$
\om_0=\pi^*\tau+d(R\Theta)=(1-R)\pi^*\tau+dR\wedge \Theta
$$
which is symplectic on $\{R<1\}\subset \cl$ and degenerates on $\{R=1 \}$ although it remains symplectic fiberwise. Let  $\spb(N,\tau):=\{R\leq 1\}_{/\sim}$ where $(z,1,\theta)\sim (z,1,\theta')$ for all $\theta,\theta'\in S^1$. As a smooth manifold, it is the projectivization of $\cl$, so we call it a {\it symplectic projective bundle}.  Since $\Theta$ is $S^1$-invariant, $\om_0$ extends smoothly to $\spb(N,\tau)$, although not to a symplectic form on the whole space: it degenerates along the "section at infinity" $N_\infty:=\{R=1\}_{/\sim}$. Nonetheless, this not really symplectic pair $(\spb(N,\tau),\om_0)$ will be convenient for our purpose. We still denote $N_0:=\{R=0\}\subset \spb(N,\tau)$ the zero-section. Notice that 
$$
\pd([\om_0])=[N_0]. 
$$
 Also,  the almost complex structure $J$ on $\sdb(N,\tau)$ defined in \S \ref{s:sdbdef} extends smoothly to an almost complex structure on $\spb(N,\tau)$, still denoted $J$. Indeed, the restriction of the almost complex structure on the fibers were defined to be the pull-back of the complex structure on the fibers of $\cl$ by the map $(R,\theta)\mapsto \big(\frac R{1-R},\theta\big)$, which extends to the complex structure on $\P^1$. Moreover, the horizontal planes $\ker \Theta$, preserved by $J$ by definition, tend to $TN_\infty$, so $J$ extends smoothly to $N_\infty$, making it a holomorphic submanifold, the map  $\pi:(N_\infty, J)\to 
 (N,j)$ being  a biholomorphism. 
 
Summarizing, given an integral symplectic manifold $(N,\tau)$, we have defined a sphere bundle $\pi_N:M:=\spb(N,\tau)\to N$, which comes along with several structures, namely: 
\begin{itemize}[leftmargin=*]
\its a closed $2$-from $\om_0$  with $[\om_0]\in H^2(M,\Z)$.
\its a smooth horizontal distribution $H^N$ on $M$ associated to a connexion form $\Theta$ defined on $M\priv (N_0\cup N_\infty)$, that coincides with the tangent  spaces of $N_0$ and $N_\infty$ along these submanifolds,
\its an almost complex structure $J$ that preserves the vertical and horizontal distributions, $N_0$ and $N_\infty$. 
\its a radial coordinate $R$ on $M$, that vanishes on the zero section $N_0$ and equals $1$ on $N_\infty$,
 \its we also have two vector fields $\frac \partial{\partial R}$ and $\frac \partial{\partial \theta}$ on $M\priv (N_0\cup N_\infty)$ and $M$, respectively, tangent to the fibers of $M$ ($\nf\partial{\partial \theta}$ vanishes along $N_0$ and $N_\infty$). 
\end{itemize}

\subsection{Line bundles over  projective bundles.}\label{s:spbb}
Let $M:=\spb(N,\tau)\overset {\pi_N}\to N$. Since $[\om_0]\in H^2(M,\Z)$, there exists a complex line bundle $\cl\overset\pi\to M$ with $c_1(\cl)=[\om_0]$ and a connexion $\nabla$ with curvature $-i\om_0$ that extends the connection $\nabla$ on $\cl\to N_0$, seen as $\cl\to N$. Indeed, starting with any connection $\nabla'$ on $M$ with curvature $-i\om_0$, $\nabla'_{|N_0}-\nabla$ is a closed $1$-form $A$ (because $\nabla'_{|N}$ and $\nabla$ have the same curvature) and $\nabla:=\nabla'-\pi_N^*A$ is another connection on the whole of $M$ with the same curvature $-i\om_0$ that coincides with $\nabla$ on $N_0$. This specific connection gives rise to connections $\nabla^k$ on $\cl^k\to M$, as defined in proposition \ref{p:nablak}, that also restrict to the product connection $\nabla^k$ on $\cl^k\to N_0\approx N$. Since there can be no ambiguity, we drop the superscript $k$ from the notation $\nabla^k$ in order to lighten the notation.     
As explained in section \ref{s:qh}, the almost complex structure $J$ provides a splitting of the connexion into its $(1,0)$ and $(0,1)$ parts, defined by 
$$
\nabla'=\frac 12(\nabla-i\nabla_{J\cdot}) \hspace{,5cm}\text{and} \hspace{,5cm} \nabla''=\frac 12(\nabla+i\nabla_{J\cdot}). 
$$  
This splitting allows to define holomorphic or quasi-holomorphic sections (those that satisfy $\nabla''s=0$ or $|\nabla''s|<\kappa|\nabla's|$, respectively). 

The aim of this paragraph is to give the basic properties of this bundle and its sections.  We start with some more or less obvious observations, that we will use several times in our computations. 

\begin{remarks}\label{rks:nabla}$ $
\begin{enumerate}[leftmargin=*]
\item[(i)]
Given any vector field $u \in \Gamma^\infty (TN)$, the vector field $U=\pi_N^*u \in \Gamma^\infty(TM)$ uniquely defined 
by $U \in H^N$ and $\pi_{N*}U=u$ commutes with $\frac\partial{\partial R}$ (in the sense of the Lie bracket)
because, since the connection defined by~$\Theta$ is linear, radial dilations preserve~$H^N$.
Thus, any vector $u \in T_pN$ can be extended to a smooth horizontal vector field on~$M$ 
that commutes with $\frac\partial{\partial R}$. As a result, 
any horizontal vector $U(p) \in T_pM$ can be extended to a smooth horizontal vector field 
that commutes with $\frac\partial{\partial R}$ and whose value at~$\pi_N(p)$ is $\pi_{N*}U(p) \in T_{\pi(p)}N$.

\item[(ii)]
Since $\Theta$ is a Hermitian connection form, the $\Theta$-parallel transport preserves $R$, so 
$$
dR(U)=0 \,\hspace{,5cm} \forall \,U\in H^N. 
$$ 
This implies in turn, taking into account the explicit formula \eqref{e:omsdb} for $\om_0$, that $H^N$
is symplectic-orthogonal to the fibers $F$ of $\pi_N$.  Also, since $J$ preserves $H_N$, 
$$
dR(JU)=0 \,\hspace{,5cm} \forall \,U\in H^N. 
$$


\item[(iii)]
For all $U\in H^N$, $\nabla_{\frac\partial{\partial R}}\nabla_U = \nabla_U\nabla_{\frac\partial{\partial R}}$.
Indeed, the difference between these two operators is the $0$-order operator $-ik\om_0(\frac\partial{\partial R},U)$, that vanishes by 
 (ii). 
 
 \item[(iv)] In particular, if $\nabla_{\frac\partial{\partial R}}s=0$ and if $U \in \Gamma^\infty(H^N)$ with $[\frac\partial{\partial R},U]=0$, 
then $\nabla_Us$ is ${\frac\partial{\partial R}}$-parallel. 

\item[(v)]
Since $J$ preserves $H^N$, (iii) implies that 
$\nabla_{\frac\partial{\partial R}}\nabla'_U=\nabla'_U\nabla_{\frac\partial{\partial R}}$ and 
$\nabla_{\frac\partial{\partial R}}\nabla''_U=\nabla''_U\nabla_{\frac\partial{\partial R}}$
whenever $U\in H^N$. And if $U\in \Gamma^\infty(H^N)$ with $\left[\frac\partial{\partial R},U\right]=0$, we get as in (iv) 
that if $\nabla_{\frac\partial{\partial R}}s=0$, both $\nabla'_Us$ and $\nabla''_U s$ are parallel, because by definition, the vector field $JU\in H^N$ also lifts the vector $j\pi_{N*}U$, and therefore also commutes with ${\frac\partial{\partial R}}$.

\item[(vi)]
$\nabla_{\frac\partial{\partial R}}\nabla_{\frac\partial{\partial \theta}}-\nabla_{\frac\partial{\partial \theta}}\nabla_{\frac\partial{\partial R}}=-ik$.

\item[(vii)]
As a result of (vi), if  $s$ is a section of $\cl$ over a \nbd of $N_0$ that verifies $\nabla_{\frac\partial{\partial R}}s=0$,
$$
\nabla_{\frac\partial{\partial \theta}}s=-ikRs.
$$
Indeed,  
$$
\nabla_{\frac\partial{\partial R}}\big(\nabla_{\frac\partial{\partial \theta}}s\big)=-iks=\nabla_{\frac\partial{\partial R}}(-ikRs)
$$
and both sections $\nabla_{\frac\partial{\partial \theta}}s$ and $ikRs$ vanish at~$R=0$ 
(recall that $\frac\partial{\partial \theta} \to 0$ as $R\to 0$). 

Similarly, if $s$ is a section defined on a \nbd of $N_\infty$ that verifies $\nabla_{\frac\partial{\partial R}}s=0$,
$$
\nabla_{\frac\partial{\partial \theta}}s=ik(1-R)s.
$$

\end{enumerate}
\end{remarks}

In the following, $P^R$ stands for the $\nabla$-parallel transport along the radial vector field $\frac \partial{\partial R}$. Notice that if a section $\sigma$ of $\cl_k$ is given along $N_0$ or $N_\infty$, the section $P^R\sigma$ provides a section of $\cl_k$ along $M\priv N_\infty$ or $M\priv N_0$, respectively. 

\begin{lemma}\label{l:secNM} Let $\sigma_0:N_0\to \cl_k$ and $\sigma_\infty:N_\infty\to \cl^k$ be sections of $\cl_k$ over $N_0$ and $N_\infty$, respectively.
The sections 
$$
s_0:=(1- R)^{\frac k2}P^R\sigma_0 \hspace*{,4cm}\text{and} \hspace*{,4cm} s_\infty:=R^{\frac k2}P^R\sigma_\infty
$$ 
are well-defined over $M$, holomorphic in restriction to the fibers, and 
$$
\begin{array}{l}
\nabla s_0=\big(-\frac k{2(1-R)}dR-ikR\Theta\big)s_0+(1-R)^{\frac k2}P^R\nabla_{\pi_{N*}\centerdot}\sigma_0\\
\nabla s_\infty=\big(\frac k{2R}dR+ik(1-R)\Theta\big)s_\infty+R^{\frac k2}P^R\nabla_{\pi_{N_\infty*}\centerdot}\sigma_\infty,\\
\end{array}
$$
where in the last formula, $\pi_{N_\infty}$ is the natural projection $M\priv N_0\to N_\infty$. 
\end{lemma}
\noindent{\it Proof:}  The computations are the same for $s_0$ and $s_\infty$, so we only deal 
 with the section $s_0$, already obviously defined over $M\priv N_\infty=M\priv\{R=1\}$. Since $P^R$ is an isometry, $P^R\sigma_0$ has bounded norm, while $1-R\to 0$ when $R\to 1$. Thus, $s_0$ extends continuously by $0$ to $N_\infty$. The extension is moreover of class $\cc^{\frac k2}$, hence $\cc^1$ because $k\geq 2$. Let us  compute $\nabla s_0$. First, since $\nabla_{\frac\partial{\partial R}}P^R\sigma_0=0$,
$$
\nabla_{\frac\partial{\partial R}}s_0=-\frac k2(1- R)^{\frac k2-1}P^R\sigma_0=-\frac k{2(1-R)}s_0.
$$
Also,  remark \ref{rks:nabla}(vii) gives
$$
\nabla_{\frac\partial{\partial \theta}}s_0=(1- R)^{\frac k2}\nabla_{\frac\partial{\partial \theta}}P^R\sigma_0=(1-R)^{\frac k2}(-ikRP^R\sigma_0)=-ik Rs_0.
$$
Now if $U\in H^N(p)$ for some point $p=(z,R,\theta)$ (in some trivialization chart), we extend it to a vector field $U\in \Gamma^\infty(H^N)$ 
that commutes with ${\frac\partial{\partial R}}$ (remark \ref{rks:nabla} (i)). In particular, it preserves the level sets of $R$, and  by remark \ref{rks:nabla} (v), we have $\nabla_U P^R\sigma=P^R\nabla_u\sigma$, where $u:=\pi_{N*}U$. So 
$$
\nabla_Us=(1- R)^{\frac k2}P^R\nabla_u\sigma.
$$ 
Altogther, these computations already prove the claimed formula for $\nabla s_0$. We finally check the holomorphicity of $s_0$ in the fibers, which amounts to checking that   $\nabla''_{\frac\partial{\partial R}}s_0=0$, because $\nabla''s_0$ is a $(0,1)$-form. Then,
$$
\nabla''_{\frac\partial{\partial R}}s_0  =\frac 12\big(\nabla_{\frac\partial{\partial R}}s_+i\nabla_{J{\frac\partial{\partial R}}}s_0\big)
  =\frac 12\big(\nabla_{\frac\partial{\partial R}}s_0+\frac i{2R(1- R)}\nabla_{\frac\partial{\partial \theta}}s_0\big)
  = 0. \square
$$

\begin{lemma}\label{l:parsecinf} There exist non-vanishing parallel sections of $\cl^k\to M$ above $N_\infty$. 
\end{lemma}
\noindent {\it Proof:} Since the connection $\nabla$ on $\cl^k$ is the natural connection on $\cl^k$ issued from the connection $\nabla$ on $\cl$, the proposition reduces to the case $k=1$ by proposition \ref{p:nablak}\,(1), which we assume henceforth. As a result, $\sdb(N,\tau)$ both lies in $M$ and in the restriction of $\cl\overset \pi\lra M$ to $N_0$. 

Since $\nabla$ has curvature $-i\om_0$ that vanishes along $N_\infty$, the $\nabla$-parallel transport has no monodromy along $N_\infty$. Though, it may have a holonomy, which is a cohomology class $\vartheta\in H^1(N_\infty,S^1)$, defined by the angle of rotation  of the $\nabla$-parallel  transport along a path representing a given homology class. The content of our claim is that $\vartheta$ vanishes. To see this, we fix a class $\alpha\in H_1(N_\infty)\overset {\pi_{N*}}\approx H_1(N_0)$, and let $\gamma:[0,1]\to N_0$ be a loop in this homology class, along whom the parallel transport in $\cl$ is the identity. Such a loop exists because the curvature is non-degenerate on $N_0$.  Let $e$ be a unitary vector of $\cl_{\gamma(0)}$ and $e_t \in \cl_{\gamma(t)}$ its parallel transport along $\gamma$. Then $(e_t)$ can be seen as a loop in $\{R=1\}\subset \sdb(N,\tau)\subset \cl$, that projects to $N_\infty$ in $\spb(N,\tau)$, giving a loop $(\wdt \gamma_t)$ along which we compute  the parallel transport. Define 
\fonction{\Phi}{[0,1)\times S^1}{\sdb(N,\tau)}{(R,t)}{Re_t.} 
The map $\Phi$ parametrizes a (half-closed) annulus in $\sdb(N,\tau)$ that  compactifies to an embedding still denoted $\Phi$ of the  closed annulus in $M=\spb(N,\tau)$, whose upper boundary $\Phi(1,t)$ is the loop $\wdt \gamma_t$. This already shows that $(\wdt \gamma_t)$ is a loop  in the class $\alpha$. By linearity of the connection, $Re_t$ is the parallel transport of the vector $Re$ along $\gamma$, so $Re_t=P^{\gamma(t)}Re\in \cl$, which implies in turn the vanishing of 
$\Phi^*\om_0$. Indeed, 
$$
\Phi^*\om_0\big(\frac\partial{\partial R},\frac\partial{\partial t}\big)=\om_0\big(\frac\partial{\partial R},R\frac d{dt}P^{\gamma(t)}e\big)=0
$$
because $\frac d{dt}P^{\gamma(t)}e\in H^N$ and by remark \ref{rks:nabla}(ii), $H_N$ is $\om_0$-orthogonal to the fibers. 

We now see $e_t$ as a parallel section of $\cl\to M$ along $\gamma$. Since the curvature of $\nabla$ vanishes on $\im \Phi$, which retracts to $\gamma$, $e$ extends to a parallel section $E$ along $\Phi$, unitary because $\nabla$ is hermitian. The restriction of $E(1,t)$ to $\Phi(\{R=1\})$ provides a non-vanishing parallel section along $\wdt \gamma_t$, showing that the parallel transport along this curve is the identity.
The holonomy  therefore vanishes on the class $\alpha$, which stood for any class in $H_1(N_\infty)$.\cqfd

\begin{figure}[h!]
\begin{center} 
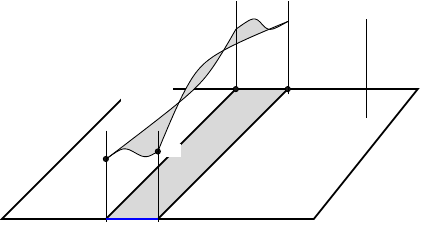
\caption{Illustration of the proof of lemma \ref{l:parsecinf}.}
\label{fig:lift}
\end{center}
\end{figure}\vspace{-,3cm}


\section{Working out an explicit section of $\pi:\cl^k\to \spb(N,\tau)$.}\label{s:sec}
Let now $\sigma:N\to \cl^k$  a quasi-holomorphic section transverse to the zero section, with $\ln|\sigma|$ Morse.
 Since $N$ embeds into $M=\spb(N,\tau)$ as the zero-section, $\sigma$ can be considered a section of $\cl^k\to N_0\subset M$.  As we have seen, we have a natural extension
$$
s_0:=(1-R)^{\frac k2}P^R\sigma:M\to \cl^k, 
$$
that obviously verifies  $Z(s_0):=s_0^{-1}(0)= \pi_N^{-1}\{\sigma=0\}\cup N_\infty$. The latter is not a smooth manifold, so $s_0$ is not transverse to the zero section. We now describe  the properties of a specific perturbation of $s_0$.

Let $\sigma_\infty:N_\infty\to \cl^k$ be a parallel section of norm $1$ provided by lemma \ref{l:parsecinf}. It is obviously $J$-holomorphic and by lemma \ref{l:secNM},  can be extended to the holomorphic section 
$s_\infty:=R^{\frac k2} P^R\sigma_\infty$ over $M$, that vanishes exactly along $N_0$. For $\eps>0$, we define 
$$
s_\eps:=s_0+\eps s_\infty. 
$$
A warning is needed here. We think of $\eps$ as a deformation parameter, hence small, but we will see in the course of the argument that this smallness is not mandatory.  And  we will need it large when  $k=2$. 
\begin{claim}\label{c:quasihol} If $\sigma$ is $\kappa$-quasi-holomorphic,  so is $s_\eps$, for the same constant $\kappa$. 
\end{claim}
\noindent{\it Proof:} Both $s_0$ and $s_\infty$ are $J$-holomorphic in restriction to the fibers and $\nabla s_\infty$ vanishes on the horizontal spaces $H^N$. Decomposing  a vector $X\in TM$ along its horizontal and vertical components $X^H$ and $X^F$, we compute:
$$
\begin{array}{llr}
|\nabla''_Xs_\eps| &=|\nabla''_{X^H} s_0|&\\
 & =|(1-R)^{\frac k2}P^R\nabla''_{\pi_{N*}X^H}\sigma| &\text{ by lemma \ref{l:secNM}}\\
 & = (1-R)^{\frac k2}|\nabla''_{\pi_{N*}X^H} \sigma | &\\
 & \leq \kappa (1-R)^{\frac k2}| \nabla'_v \sigma| & \text{ for some } v \text{ with }|v|= |\pi_{N*}X^H|\\
 & = \kappa |(1-R)^{\frac k2} P^R\nabla'_{V}\sigma| & \text{ with } V:=\pi_N^*v\in H^N, \\
 & = \kappa |\nabla'_V s_0| & \\
  & =\kappa| \nabla'_V s_\eps| & \text{ because } \nabla s_\infty \equiv 0 \text{ on }H^N.   
\end{array} 
$$
Moreover, by definition of the metric, 
$$
|V|=|v|= |\pi_{N*}X^H|=|X^H|\leq |X|. 
$$
Thus for all $X$ in $T_pM$, we have found a vector $V\in T_pM$ with smaller norm, such that $|\nabla''_Xs_\eps|\leq \kappa |\nabla'_Vs_\eps|$. This precisely means that $|\nabla''s_\eps|\leq \kappa |\nabla' s_\eps|$, so $s_\eps$ is indeed $\kappa$-quasi-holomorphic. \cqfd

\begin{figure}[h!]
\begin{center} 
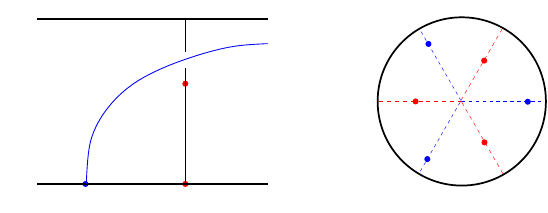
\caption{Zero set of $s_\eps$, and critical points of $\ln|s_\eps|$ when $k\geq 3$.}
\label{fig:Zcrit}
\end{center}
\end{figure}

We now investigate both the zero set and the critical set of $|s_\eps|^2$, called $Z(s_\eps)$ and $\crit(|s_\eps|^2)$, respectively. Figure \ref{fig:Zcrit} p.\,\pageref{fig:Zcrit} summarizes the results obtained in the remaining of this section. One checks immediately the following:
$$
\begin{array}{l}
Z(s_\eps)\cap N_0=Z(\sigma):=\{\sigma=0\}\subset N_0, \hspace{,5cm} Z(s_\eps)\cap N_\infty=\emptyset,\\
\crit |s_\eps|^2\supset Z(s_\eps),  \hspace{,5cm}  \crit |s_\eps|^2\cap \pi_N^{-1}(Z(\sigma))=Z(\sigma). 
 \end{array}
$$
We first focus  on the trace of these subsets in $M\priv (N_0\cup N_\infty)=\sdb(N)\priv N_0$. 
Notice first that since the parallel transport is a unitary map,  $P^R\sigma_\infty$ has constant norm $1$, and we can write 
$$
P^R\sigma=|\sigma|e^{i\psi} P^R\sigma_\infty 
$$
where $\psi$ is a well-defined $S^1(2\pi)$-valued function on $M\priv (N_0\cup N_\infty\cup Z(s_0))$.

\begin{claim}\label{cl:psi} 
The function $\psi$ verifies:
\begin{itemize}
\its $\frac{\partial \psi}{\partial R}=0$ and $\frac{\partial \psi}{\partial \theta}=-k$,
\its The vectors $(\pi_{N*}\nabla \psi,\nabla \ln|\sigma|)$ form a free family as soon as one of them do not vanish, and the vanishing of one of these vectors implies the vanishing of the other.
\end{itemize}
\end{claim}
\noindent{\it Proof:} Let $(z,R,\theta)$ be a trivialization chart of $\sdb(N,\tau)\to N$ giving coordinates on a part of $\sdb(N,\tau)$. By definition, we have 
\begin{equation}\label{eq:clpsi1}
s_0=\left(\frac{1-R}{R}\right)^{\frac k2}|\sigma| e^{i\psi}s_\infty=:\left(\frac{1-R}{R}\right)^{\frac k2}e^{-ik\theta}|\sigma| e^{i\psi'}s_\infty. 
\end{equation}
One checks without difficulty that in each fiber $w:=\left(\frac{R}{1-R}\right)^{\frac 12} e^{i\theta}$ is  holomorphic. Since $s_0$ and $s_\infty$ are fiberwise holomorphic, the restriction of  $|\sigma(z)|e^{i\psi'(z,R,\theta)}$ to each fiber is holomorphic as well, with values in a circle. This function is 
therefore fiberwise constant, so $\psi'=\psi'(z)$, and our first bullet is proved.

Let now $U\in H^N$. By \eqref{eq:clpsi1}, since $\nabla_Us_\infty=0$, 
$$
\nabla_Us_0=\left(\frac{1-R}R\right)^{\frac k2} d(\sigma e^{i\psi})(U)s_\infty=(d\ln|\sigma|(U)+id\psi(U))s_0. 
$$
Since $|\nabla''s_{0|H^N}|\leq \kappa |\nabla' s_{0|H^N}|$, lemma \ref{l:qhsec} shows the second assertion. \cqfd

Let us introduce the function \fonction{x}{M\priv N_0}{\R^+}{(z,R,\theta)}{\ds \left(\frac{1-R}R\right)^{\frac k2} \frac{|\sigma(z)|}{\eps}.}
Then,
$$
\begin{array}{ccl}
s_\eps &=& \eps(1+xe^{i\psi})s_\infty,\\
|s_\eps|^2 &=
 & \eps^2 R^{k}\left| 1+xe^{i\psi} \right|^2
   = \eps^2R^{k} (x^2+2x\cos\psi +1),
\end{array}
$$ 
so that 
$$
\begin{array}{lr}
d|s_\eps|^2 = & k\eps^2R^{k-1}(x^2+2\cos \psi x +1)dR\\ & +2\eps^2R^k(x+\cos \psi) dx \\ & -2\eps^2R^kx\sin \psi d\psi. 
\end{array}
$$
Further, 
\begin{equation}\label{eq:dx}
dx=d\left[\left(\frac {1-R}{R}\right)^{\frac k2}\frac{|\sigma|}\eps\right]=-\frac{kx}{2R(1-R)}dR+xd\ln|\sigma|. 
\end{equation}
After a straightforward computation, 
we get on $M\priv N_0$:
\begin{equation}\label{eq:grads}
\begin{array}{lr}
d |s_\eps|^2 &  = \frac{k\eps^2R^{k-1}}{1-R}\big(-Rx^2+(1-2R)\cos \psi x+1-R\big)dR\\
 & +2\eps^2R^{k}(x+\cos \psi)x\, d \ln |\sigma|\\
  & -2\eps^2R^kx\sin \psi\, d \psi
\end{array}
\end{equation}
and 
\begin{equation}\label{eq:|s|2}
|s_\eps|^2  = \eps^2R^{k} (x^2+2x\cos\psi+1). 
\end{equation}
The latter already gives the vanishing set of $s_\eps$:
$$
\begin{array}{ll}
Z(s_\eps)&:= \{|s_\eps|=0\}=\{R=0,\; \sigma=0\}\cup\{\cos\psi=-1,\; x=1\}\\
 & \ds =\{\theta=\frac{\pi-\psi'(z)}k\,\left[\frac {2\pi}k\right],\; R=\big(1+\left(\frac \eps{|\sigma(z)|}\right)^{\frac 2k}\big)^{-1}\}\cup Z(\sigma) 
 \end{array}
$$
where we recall that $Z(\sigma):=\{\sigma=0\}\subset N_0$. 
In each fiber for which $\sigma(z)\neq 0$, there are exactly $k$ zeroes of $s_\eps$ which are  equi-distributed on a circle of fixed radius. This radius gets closer to $1$ as $\eps\to 0$ and closer to $0$ when $\eps\to \infty$, uniformly on compact sets of $\{\sigma(z)\neq 0\}$. On the other hand, 
$$
Z(s_\eps)\cap \pi_N^{-1}(Z(\sigma))=\{\sigma(z)=0,R=0\}=Z(\sigma)\subset N_0.
$$
\begin{lemma}\label{l:trvanish} The section $s_\eps$ vanishes transversally in $\sdb(N,\tau)\priv N_0$. 
\end{lemma}
\noindent{\it Proof:} We have seen that $Z(s_\eps)\priv N_0$ lies in the complement of $Z(s_0)$, where  $s_\eps=\eps(1+xe^{i\psi})s_\infty$. Since $\nabla_{\frac\partial{\partial R}}s_\infty=\frac k{2R}s_\infty$, equation \eqref{eq:dx} gives: 
$$
\begin{array}{ll}
\nabla_{\frac\partial{\partial R}}s_\eps&=\eps\big(\frac k{2R}(1+xe^{i\psi})+e^{i\psi}\frac {\partial x}{\partial R}\big)s_\infty\\ & =\eps\big(\frac k{2R}(1+xe^{i\psi})-\frac {kxe^{i\psi}}{2R(1-R)}\big)s_\infty\\ &=\frac{\eps k}{2R(1-R)}(1-R-Rxe^{i\psi})s_\infty.  
\end{array}
$$
At a point of $Z(s_\eps)\priv N_0$, $x=1$ and $e^{i\psi}=-1$ so 
$$
\nabla_{\frac\partial{\partial R}}s_\eps=\frac{\eps k}{2R(1-R)}s_\infty\neq 0.
$$
As a result, the radial covariant derivative of $s_\eps$ at such a point does not vanish. Since $s_\eps$ is fiberwise holomorphic, this implies that $\nabla s_\eps$ has rank $2$ when restricted to the fiber directions at such a point, so $\nabla s_\eps$ is surjective. \cqfd

Using formula \eqref{eq:grads} we now investigate $\crit |s_\eps|^2\priv Z(s_\eps)$ in $M\priv (N_0\cup N_\infty)$. Let  $p\in \crit |s_\eps|^2\priv Z(s_\eps)$, $z:=\pi_N(p)$ and $R:=R(p)\notin\{0,1\}$. We have already noticed that $p\notin Z(s_0)=\pi_N^{-1}(Z(\sigma))$ so $x,\psi$ are well-defined around $p$. 
By claim \ref{cl:psi}, the $1$-forms $(dR,d\psi)$ form a free family and $d\ln|\sigma| \notin \span(dR,d\psi)$ at $p$ whenever $d\ln|\sigma|(p)\neq 0$.  The vanishing of the component of $d|s_\eps|^2$ along $d\psi$ at $p$ in equation \eqref{eq:grads} therefore gives
$$
\sin\psi=0,
$$  
which in turn implies $\cos \psi=\pm 1$. Since $R(p)\notin \{0,1\}$, the vanishing of the radial component of $d|s_\eps|^2$ implies that $x\neq 0$. Thus, the vanishing of the last component  implies that 
either $x+\cos \psi=0$, or $d \ln|\sigma|=0$. Since $\cos\psi=\pm 1$ and $x\geq 0$, the first case is not compatible with our assumption that $p\notin Z(s_\eps)$, since it implies that $\cos \psi=-1$ and $x=1$. So we must have $d \ln|\sigma|=0$. Finally, the vanishing of the radial component gives 
$$
Rx^2-(1-2R)\cos \psi\, x -(1-R)=0, 
$$
which gives two solution compatible with $x\geq 0$:
$$
\begin{array}{ll}
\cos \psi=1,&  x=\frac{1-R}{R}\\
\cos \psi=-1 & x=1.
\end{array}
$$
The second line describes the set $Z(s_\eps)$. The first line  describes a genuine critical point when $k\geq 3$, but when $k=2$ it amounts to 
$$
\frac{1-R}{R}\frac{|\sigma|}\eps =\frac {1-R}R,
$$
which has no solution under the condition $\eps>\max |\sigma|$,  henceforth assumed. 
In summary, under our conditions on $(k,\eps)$:
$$
\begin{array}{l}
 \crit |s_\eps|^2\cap M\priv (N_0\cup N_\infty)=\\
\hspace{1,5cm} Z(s_\eps)\cup\big\{(z,R,\theta)\; |\;\left\{\begin{array}{l} d\ln|\sigma|(z)=0, \\ \cos(k\theta+\psi'(z))=1\\  R=\big(1+\left(\frac{\eps}{|\sigma(z)|}\right)^{\frac 2{k-2}}\big)^{-1}\end{array}\right.\big\} \text{ when }k\geq 3, \\
\hspace{1,5cm}Z(s_\eps) \text{ when } k=2.
 \end{array}
$$
In other terms, when $k\geq 3$ there are exactly $k$ critical points for $|s_\eps|^2$ on each fiber of the critical points of $|\sigma|$, equi-distributed along a circle. These fibers also intersect the vanishing set of $s_\eps$ along $k$ points located on the same diameters as the critical points, only in the other direction from the origin. 

\vspace{-,3cm}

\begin{lemma}\label{l:inds3} 
When $k\geq 3$, the function $-|s_\eps|^2$ is Morse at any  $p_0\in \crit(|s_\eps|^2)\priv  \big(Z(s_\eps)\cup N_0\cup N_\infty\big)$ and its Morse index is one more than the index of $\sigma$ at the point $\pi(p_0)$. 
\end{lemma}
\noindent{\it Proof:} Recall the formula \eqref{eq:grads} for the differential of $|s_\eps|^2$:
$$
\begin{array}{ll}
-d|s_\eps|^2=&\frac{k\eps^2R^{k-1}}{1-R}\big(Rx^2-(1-2R)x\cos\psi-1+R\big)dR\\
 & -2\eps^2R^k(x+\cos\psi)xd\ln|\sigma|\\
 & -2\eps^2R^kxd \cos\psi. 
\end{array}
$$
At a critical point $p_0\in \sdb(N,\tau)\priv Z(s_\eps)$ that projects to $z_0$, we have $\cos \psi(p_0)=1$, $d\ln|\sigma|(z_0)=0$, $x=\frac{1-R}R$. A lengthy but straightforward computation gives:
$$
\begin{array}{ll}
-d^2|s_\eps|^2=&-\frac{k(k-2)\eps^2R^{k-3}}{2(1-R)}dR^2\\
 & +\eps^2R^{k-1}(1-R)d\psi^2\\ 
 & -2\eps^2R^{k-2}(1-R)d^2\ln|\sigma|.
 \end{array}
$$
By claim \ref{cl:psi}, $T^*_{p_0}M=\la dR\ra\oplus\la d\psi(p_0)\ra\oplus H_N^*(p_0)$, so this decomposition shows that $-|s_\eps|^2$  is Morse at $p_0$ with index   one more than that of $-\ln|\sigma|$ at $z_0$.\cqfd

\begin{lemma}\label{l:inds2} When $k=2$, $\crit |s_\eps|^2\cap N_0=\crit |\sigma|^2$. At a point $p\in \crit|\sigma|^2\priv Z(\sigma)\subset N_0$, $-|s_\eps|^2$ is Morse with index one more than the index of $-|\sigma|^2$. 
\end{lemma}
\noindent{\it Proof:} When $k=2$, 
$$
s_\eps=(1-R)P^R\sigma+\eps RP^R\sigma_\infty.  
$$
Fix a point $z_0\in N$, coordinates $(z,w)$ on $\pi_N^{-1}(\op(z_0))\priv N_\infty$, $w=re^{i\theta}=x+iy$, $R=r^2$. By claim \ref{cl:psi},
$$
\begin{array}{ll} 
|s_\eps|^2(z,w) &=(1-R)^2|\sigma(z)|^2+\eps^2R^2+2\eps R(1-R)|\sigma(z)|\cos (2\theta +\psi'(z))\\ &=|\sigma(z)|^2+O(R). \end{array}
$$
This already shows that $\crit |s_\eps|^2\cap N_0=\crit|\sigma|^2$. 
Let now $z_0\in \crit|\sigma|^2$, write $(z,w)=(z_0+h,w)$ and modify the trivialization in order to make $\psi'(z_0)=0$. We then have 
$$
\begin{array}{rl}
|s_\eps|^2(z,w)=&  |\sigma(z)|^2-2R|\sigma(z)|^2+2\eps|\sigma(z)|R\re(e^{i2\theta+\psi'(z)})+O(R^2)\\
= &   |\sigma(z)|^2-2R|\sigma(z_0)|^2+2\eps|\sigma(z_0)|\re(Re^{i2\theta})+o(|(h,w)|^2)\\
=&  |\sigma(z)|^2-2R|\sigma(z_0)|^2+2\eps|\sigma(z_0)|\re(w^2)+o(|(h,w)|^2\\
=& |\sigma(z)|^2-2|\sigma(z_0)|^2(x^2+y^2)+2\eps|\sigma(z_0)|(x^2-y^2)+o(|(h,w)|^2\\
=&  |\sigma(z)|^2+2\eps|\sigma(z_0)|\left((1-\frac{|\sigma(z_0)|}{\eps})x^2-(1+\frac{|\sigma(z_0)|}{\eps})y^2\right)+o(|(h,w)|^2).\\
\end{array}
$$
Since $\eps>\max |\sigma|$, this order $2$ expansion of $|s_\eps|^2$ at $(z_0,0)$ gives the desired result, and also shows that the stable direction in the fiber is $\{y=0\}$ , which corresponds to $\cos \psi=1$. \cqfd

Although we won't formally need it, notice that it is easy to see that $\crit |s_\eps|^2\cap (N_0\cup N_\infty)=\crit |\sigma|^2\cup N_\infty$ and that $|s_\eps|^2$ is never Morse on $N_\infty$ and Morse on $N_0$ when $k=2$ but not when $k\geq 3$. This could be an issue when studying the Liouville flow associated to the section $s_\eps$ (in the next section), if we did not have the following remarks on the variations of $|s_\eps|^2$ near $N_0\cup N_\infty$. 

\begin{lemma}\label{l:varsinf}
The function $-|s_\eps|^2$ increases along the radial directions emanating from $N_\infty$. When $k\geq 3$, it also increases along the radial directions emanating from $N_0$. 
\end{lemma}
\noindent{\it Proof:} By formula \eqref{eq:grads}, on $R\neq 0,1$:
$$
\begin{array}{ll}
-\frac{\partial |s_\eps|^2}{\partial R} &=\alpha\big(Rx^2-(1-2R)x\cos \psi-(1-R)\big)\\  &=\beta\big(R-(1-2R)x^{-1}\cos \psi-(1-R)x^{-2}\big)
\end{array}
$$
where $\alpha$ and $\beta$ are positive functions. Recall also that 
$$
x=\left(\frac {1-R}R\right)\frac{|\sigma|}\eps,
$$
so $x\underset{R\simeq 1}\sim (1-R)^{\nf k2}$ and $x^{-1}\underset{R\simeq 1}\sim R^{\nf k2}$. Thus when $k\geq 3$,
$$
Rx^2-(1-2R)x\cos \psi-(1-R)\underset{R\simeq 1}\sim -(1-R)<0
$$
and 
$$
R-(1-2R)x^{-1}\cos \psi-(1-R)x^{-2}\underset{R\simeq 0}\sim R>0,
$$
so indeed $-|s_\eps|^2$ increases along the radial directions emanating from $N_0$ and $N_\infty$. When $k=2$, 
$$
Rx^2-(1-2R)x\cos \psi-(1-R)\underset{R\simeq 1}= -(1-R)\big(1 +\frac {|\sigma|}\eps\cos\psi\big))+o(1-R)
$$
and the right hand side is clearly negative under our assumption that $\eps>\max |\sigma|$. \cqfd


\section{A Liouville vector field on $\sdb(N,\tau)\priv Z(s_\eps)$ and its skeleton.} \label{s:lvf}
\begin{figure}[h!]
\begin{center} 
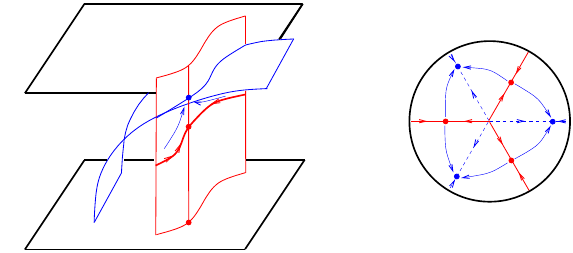
\caption{Skeleton of $s_\eps$.}
\label{fig:skeleps}
\end{center}
\begin{minipage}{0,8\textwidth}{\footnotesize The points in $\pi_N^{-1}(\skel(\sigma))$ are attracted by the fibers over $\crit(\sigma)$, but only a minority are attracted by a point of $\crit(s_\eps)$, the rest abute in finite time at $Z(s_\eps)\cap \pi_N^{-1}(\skel(\sigma))$.}\end{minipage}
\end{figure}

We keep the same framework and notation as in the previous paragraph where we have constructed a section $s_\eps$ of $\cl^k\to M=\spb(N,\tau)$ out of a section $\sigma$ of the restriction of this bundle to $N_0$, and studied its vanishing set and critical points. 
As in \cite{giroux1}, define $\lambda_\eps$ on $M\priv Z(s_\eps)$ by 
$$
\nabla(f \frac {s_\eps}{|s_\eps|})=(df-ik\lambda_\eps f)\frac {s_\eps}{|s_\eps|},
$$
or equivalently by 
$$
\nabla s_\eps=(d\ln|s_\eps|-ik\lambda_\eps)s_\eps=(-d\phi_\eps-ik\lambda_\eps)s_\eps,
$$
where we have put
$$
\phi_\eps:=-\ln|s_\eps|.
$$
Since $\nabla$ is a smooth connection on $M$, $\lambda_\eps$ is smooth on $M\priv Z(s_\eps)$, and since $\nabla$ has curvature $-ik\om_0$, $d\lambda_\eps=\om_0$. On $M\priv (Z(s_\eps)\cup N_\infty)=\sdb(N,\tau)\priv Z(s_\eps)$, on which $\om_0$ is non-degenerate, we have a well-defined associated Liouville vector field $X_\eps$. Notice that when we restrict $s_\eps$ to $N_0$, we get 
$$
\nabla \sigma=(d\ln|\sigma|-ik\lambda_\eps)\sigma,
$$
so $\lambda:=\lambda_{\eps|N_0}$ is the Liouville form naturally associated to the section $\sigma$. This paragraph  aims at describing the basin of attraction of $Z(s_\eps)$ under the flow of $X_\eps$, and its results are pictured on figure \ref{fig:skeleps} p.\,\pageref{fig:skeleps}. 

Since $s_\eps$ is quasi-holomorphic, lemma \ref{l:qhsec} shows that $X_\eps$ is gradient-like for $\phi_\eps$ (meaning that $X_\eps\cdot \phi_\eps\geq 0$, with equality only at critical points of $\phi_\eps$, where $X_\eps$ vanishes). By standard arguments on gradient flows, the basin of attraction of $Z(s_\eps)$ is therefore the complement of the set of points attracted by $\crit \phi_\eps$ (under $X_\eps$). By lemma \ref{l:varsinf} we already know that $N_\infty$, on which $X_\eps$ is not well-defined, is repulsive so attracts nothing. When $k\geq 3$, $\crit \phi_\eps\cap N_0$ attracts only $\skel\,\sigma\subset N_0$, and the critical points in $M\priv (N_0\cup N_\infty)$ are Morse. When $k=2$, $\phi_\eps$ has no critical point in $M\priv (N_0\cup N_\infty)$ and the critical points of $\phi_\eps$ on $N_0$ are Morse. What therefore remains to do is to identify the stable manifolds of the critical points of $\phi_\eps$ in $M\priv(N_0\cup N_\infty)$ when $k\geq 3$ and in $N_0$ when $k=2$.

\begin{lemma}\label{l:psilambda} 
$$
\pi_N^*\lambda=-\frac 1k d\psi-\Theta.  
$$
As a result, the function $\psi$ remains constant along the trajectories of $\pi_N^*X_\lambda$ (the horizontal vector field that projects to $X_\lambda$). 
\end{lemma}
\noindent{\it Proof:} We recall that $\lambda:=\lambda_{\eps|N_0}$ verifies $\nabla \sigma=(d\ln|\sigma|-ik\lambda)\sigma$. 
On one hand, lemma \ref{l:secNM} gives 
$$
\begin{array}{ll}
\nabla s_0&=\big(-\frac k{2(1-R)} dR-ikR\Theta)s_0+(1-R)^{\frac k2}P^R\nabla_{\pi_{N*\centerdot}} \sigma\\
& =\big(-\frac k{2(1-R)} dR-ikR\Theta)s_0+(1-R)^{\frac k2}P^R(\pi_N^*d\ln|\sigma|-ik\pi_N^*\lambda)\sigma\\
& =\big(-\frac k{2(1-R)} dR-ikR\Theta)s_0+(d\ln|\sigma\circ \pi_N|-ik\pi^*_N\lambda)s_0\\
\end{array}
$$
so
$$
-\im\frac{\nabla s_0}{k s_0}=\pi^*_N\lambda+R\Theta.
$$
On the other hand, writing $s_0=\eps xe^{i\psi}s_\infty$ and using 
 lemma \ref{l:secNM} again (taking into account that $\nabla\sigma_\infty=0$), we get:
$$
\frac{\nabla s_0}{s_0}=d\ln x+id\psi+\frac{\nabla s_\infty}{s_\infty}
 = d\ln x+\frac k{2R}dR+i(d\psi +k(1-R)\Theta),
$$
so 
$$
-\im \frac {\nabla s_0}{k s_0}=-\frac{d\psi}k-(1-R)\Theta. 
$$
Comparing these two expressions, we deduce that 
$\pi_N^*\lambda=-\frac 1k d\psi-\Theta$.  As a result, 
$$
0=\lambda(X_\lambda)=\pi_N^*\lambda(\pi_N^*X_\lambda)=-\frac1 k d\psi(\pi_N^*X_\lambda)-\Theta(\pi_N^*X_\lambda)=\frac1 k d\psi(\pi_N^*X_\lambda).\hspace{,3cm}\square
$$

\begin{lemma}\label{l:lambdaeps} 
$$
\lambda_\eps=\frac{x^2+x\cos \psi}{|1+xe^{i\psi}|^2} \pi_N^*\lambda+\frac{Rx^2+(2R-1)x\cos\psi-1+R}{|1+xe^{i\psi}|^2}\Theta-\frac{\sin\psi}{k|1+xe^{i\psi}|^2}dx
$$
\end{lemma}
\noindent{\it Proof:}  Using that for any section $s$ of $\cl^k$, $\re \frac {\nabla s}s=\ln|s|$ (lemma \ref{l:qhsec}) and lemma \ref{l:secNM}, we get:
$$
\begin{array}{ll}
\nabla s_\eps&=\nabla s_0+\eps\nabla s_\infty\\
 & = \big(d\ln|s_0|-ik(\pi_N^*\lambda+R\Theta)\big)s_0+\eps\big(d\ln|s_\infty|+ik(1-R)\Theta\big)s_\infty. 
\end{array}
$$
Since $s_0=\eps xe^{i\psi} s_\infty$ and $s_\eps=s_0+\eps s_\infty$, 
$$
s_\infty=\frac {s_\eps}{\eps(1+xe^{i\psi})}, \hspace{,5cm} s_0=\frac{xe^{i\psi}}{1+xe^{i\psi}}s_\eps. 
$$
Thus, 
$$
\begin{array}{ll}
 |1+xe^{i\psi}|^2\frac{\nabla s_\eps}{s_\eps}=& \big(d\ln|s_0|-ik(\pi_N^*\lambda+R\Theta)\big)(x^2+xe^{i\psi})+\big(d\ln|s_\infty|+\\ & +ik(1-R)\Theta\big)(1+xe^{-i\psi})
 \end{array}
$$
and 
$$
 \begin{array}{rl}
 |1+xe^{i\psi}|^2\im \frac{\nabla s_\eps}{s_\eps}=&x\sin \psi d\ln\left(\frac{|s_0|}{|s_\infty|}\right)-k(x^2+x\cos \psi)(R\Theta+\pi_N^*\lambda)+\\ & +k(1-R)(1+x\cos \psi)\Theta\\
 =& x\sin \psi d\ln|x|-k\big(Rx^2+(1-2R)x\cos\psi+1-R\big)\Theta-\\ & -k(x^2+x\cos\psi)\pi_N^*\lambda. 
 \end{array}
$$
This concludes the proof, since $\lambda_\eps=-\im \frac{\nabla s_\eps}{ks_\eps}$.\cqfd 

Since $X_{\pi_N^*\lambda}=\pi_N^*X_\lambda$ and $X_\Theta=\frac\partial{\partial R}$, the previous lemma gives the following explicit formula for the Liouville vector field associated to $\lambda_\eps$:
\begin{equation}\label{eq:Xeps}
\textstyle X_\eps=\frac{x^2+x\cos \psi}{|1+xe^{i\psi}|^2} \pi_N^*X_\lambda+\frac{Rx^2+(2R-1)x\cos\psi-1+R}{|1+xe^{i\psi}|^2}\frac\partial{\partial R}-\frac{\sin\psi}{k|1+xe^{i\psi}|^2}X_x,
\end{equation}
where $X_x$ denotes as usual the Hamiltonian vector field associated to the function $x$. We are now in position to compute the stable manifolds of the critical points of $\phi_\eps$ under $X_\eps$. 

\begin{corollary}\label{c:stseps} The stable manifold of $X_\eps$ at a critical point $p=(z_0,R_0,\theta_0)\in \sdb(N,\tau)\priv Z(s_\eps)$ of $\phi_\eps$ is given by:
$$
W^s\big(X_\eps,p\big)=\{(z,R,\theta)\; |\; z\in W^s(X_\lambda,z_0), \cos \psi(z,\theta)=1\}. 
$$
 \end{corollary}
\noindent{\it Proof:} 
Let $p_0:=(z_0,R_0,\theta_0)$ be a critical point of $|s_\eps|^2$  in $\sdb(N,\tau)\priv Z(s_\eps)$. We recall that it verifies:
\begin{itemize}
\its $z_0\in \crit\ln|\sigma|$, $e^{i\psi}=1$ and $x=\frac {1-R_0}{R_0}$  when $k\geq 3$,
\its $z_0\in \crit \ln |\sigma|$, $R_0=0$ and $\theta_0$ is undefined when $k=2$. 
\end{itemize}
 In both cases, let 
$$
W:=\{(z,R,\theta)\;|\; z\in W^s(X_\lambda,z_0),\cos\psi(z,\theta)=1\}. 
$$
We claim that $W\subset W^s(X_\eps,p_0)$. Let indeed $p=(z,R,\theta)\in W$ and define $p(t):=(z(t),R(t),\theta(t))$ by 
$$
\left\{\begin{array}{l}
(z(0),R(0),\theta(0))=(z,R,\theta),\\
\dot z(t)=\frac{x^2+x}{(1+x)^2}X_\lambda(z(t))=\frac x{1+x}X_\lambda(z(t))\\
\dot R(t)=\frac{Rx^2+(2R-1)x-1+R}{(1+x)^2}=\frac{Rx-1+R}{1+x}\\
\dot \psi(t)=0,
 \end{array}\right.
$$ 
for $t\geq 0$, as long as these equations make sense. Notice that the second and third lines do not involve the coordinate $\theta$, so define a differential equation on $(R,z)$. Once it is solved, $z(t)$ and $R(t)$ are fixed, and the last line, that lemma \ref{l:psilambda} allows to write 
$$
-k\dot\theta(t)+\lambda(\dot z(t))=0,
$$
fixes the derivative of $\theta(t)$, hence $\theta$ itself, and forces $\cos \psi(p(t))=1$ along the trajectory. Notice also that this last line could equivalently be written $\cos\psi(p(t))=1$, making no reference to the coordinate $\theta$, so this differential equation is well-defined on the whole of $M$, not only on a trivialization chart of $\spb(N,\tau)\to N$. Now $p(t)$ is obviously a solution of the differential equation $\dot p(t)=X_\eps(p(t))$ by formula \eqref{eq:Xeps} and  $\cos\psi(z(t),\theta(t))=1$. Moreover,
\begin{itemize}[leftmargin=*]
\its when $k\geq 3$,  $R$ increases or decreases, depending on whether $x<\frac {1-R}R$ or $x>\frac{1-R}R$, respectively, so $R$ remains in a compact region of $(0,1)$ when $t$ grows. Also, $z(t)\in W^s(X_\lambda,z_0)$, so does not approach $Z(\sigma)$, and  $\frac x{1+x}$ remains in a compact region of $\R^+\priv\{0\}$. This implies in turn that the trajectory is well-defined for all $t\geq 0$ and that $z(t)$ converges to $z_0$. As a result $p(t)$ converges to $p_0$ when $t$ goes to $+\infty$. 

\its When $k=2$, the equations on $(z,R)$ can be rewritten
$$
\left\{
\begin{array}{l}
\dot R(t)=-\frac{R}{(1-R)(1+x^{-1})}(\frac \eps{|\sigma|}-1)<0,\\
\dot z(t)=\frac 1{1+x^{-1}}X_\lambda(z(t)). 
\end{array}
\right.
$$
Thus, $R$ decreases with $t$ and since $z(t)\in W^s(X_\lambda,z_0)$, $z(t)$ does not approach $Z(\sigma)$.  Thus, 
$$
x^{-1}=\frac R{1-R}\frac{\eps}{|\sigma|}
$$
remains bounded. The trajectory is then defined for all $t\geq 0$, $z(t)$ converges  to $z_0$ and $R(t)$ converges to $0$. Thus, again $p(t)$ converges to $p_0$ when $t$ goes to $+\infty$. 
\end{itemize}
 Thus, in both cases, $p\in W^s(X_\eps,p_0)$, hence the claimed inclusion $W\subset W^s(X_\eps,p_0)$. Now it is clear that $W$ is a submanifold  of dimension one more than the index of $-\ln|\sigma|$ at $z_0$, which is the Morse index of $\phi_\eps$ at $p_0$ by lemma \ref{l:inds3} and \ref{l:inds2}. Since  $X_\eps$ is gradient-like for $\phi_\eps$, $W$ is the stable manifold of $X_\eps$ at $p_0$. \cqfd

We define now $\cb(Z(s_\eps),X_\eps)\subset \sdb(N,\tau)$ as the set of points of $\sdb(N,\tau)$ attracted by $Z(s_\eps)$ under the flow of $X_\eps$ and $\skel(s_\eps)$ its complement in $\sdb(N,\tau)$. As already pointed out, $\skel(s_\eps)$ is the set of points whose trajectories under $X_\eps$ are defined for positive time and converge to a critical point of $\phi_\eps$. 
\begin{lemma}\label{l:idskel}
Let $\Gamma_k\subset N$ be the isotropic skeleton associated to the section $\sigma:N\to \cl^k$, namely
 $$
\Gamma_k:=\bigcup_{z\in \crit |\sigma|^2\priv Z(\sigma)} W^s(X_\lambda,z).
$$ 
Then there exists a $k$-fold Legendrian lift $\Lambda_k$ of $\Gamma_k$ in the prequantization $\cp(N,\tau)$ such that 
$$
\skel(s_\eps)=\overline{\Phi^{-1}(\Lambda_k\times (0,1))},
$$
where $\Phi$ is the identification between $\sdb(N,\tau)\priv N_0$ and $S^{<0}\cp$ provided by lemma \ref{l:sdbsp} and the closure is meant in $\sdb(N,\tau)$.  
\end{lemma}
\noindent{\it Proof:} Decompose 
$$
\skel(s_\eps)=\big(\skel(s_\eps)\priv N_0\big)\cup \big(\skel(s_0)\cap N_0\big). 
$$
Corollary \ref{c:stseps} gives 
$$
\skel(s_\eps)\priv N_0  =\left\{p\in M\; |\; \left\{\begin{array}{l}  \pi_N(p)\in \bigcup_{z\in \crit |\sigma|^2\priv Z(\sigma)} W^s(X,z),\\ \cos\psi(p)=0\end{array}\right. \right\}. 
$$
 By \cite{giroux2}, this skeleton is known to be isotropic in $\sdb(N,\tau)$ because it is the skeleton of the quasi-holomorphic section $s_\eps$. The  formula on the right shows that it is also radially invariant and that it projects to $\Gamma_k$. Finally, since $\frac{\partial\psi}{\partial \theta}(p)=-k$, the fiber above a point in $\Gamma_k$ is a union of $k$-equidistributed rays. The uniqueness part of lemma \ref{l:skelift} therefore guarantees that 
$$
\skel(s_\eps)\priv N_0=\Phi^{-1}(\Lambda_k\times (0,1)),
$$  
where $\Lambda_k$ is a $k$-fold isotropic cover of $\Gamma_k$ in $\cp(N,\tau)$. On the other hand, $X_\eps$ is tangent to $N_0$, where it coincides with the Liouville flow associated to the quasi-holomorphic section $\sigma$. Thus, 
$$
\skel(s_\eps)\cap N_0=\skel(\sigma)=\Gamma_k. 
$$
It is clear that the union of these two subsets is the closure of the former in $\sdb(N,\tau)$. 
\cqfd 

In summary, we have obtained the following version of Biran decomposition theorem, whatever $k\geq 2$:
$$
\sdb(N,\tau)=\cb(Z(s_\eps),X_\eps)\sqcup \overline{\Phi^{-1}(\Lambda_k\times (0,1))}. 
$$
This decomposition is illustrated in the following figure. 

\vspace{-,3cm}

\section{Taming the Liouville form.}\label{s:taming}

We finally need to  describe a symplectic model of $\cb(Z(s_\eps),X_\eps)$. Let us introduce the notion of tameness of Liouville forms, as proposed in \cite{opshtein3,opsc}:

\begin{definition} Let $(M^{2n},\om)$ be a symplectic manifold, $\Sigma^{2n-2}\subset M$ a smooth symplectic divisor locally  Poincaré-dual to  a multiple  of the symplectic class (recall that this means that the first Chern class of its normal symplectic bundle is $\ell[\om_{|\Sigma}]$). Use the symplectic \nbd theorem to present a \nbd of $\Sigma$ in $M$ as a \nbd of the zero section in $\sdb(\Sigma,\om_{|\Sigma},\ell)$, on which we have the connection form $\Theta$ (see section \ref{s:sdbdef}). We say that a $1$-form $\lambda$ defined on $M\priv \Sigma$ is tame along $\Sigma$ if there exists a constant $a\in \R$ and a bounded $1$-form $\nu$ in $M\priv \Sigma$ such that  $\lambda=a\Theta+\nu$. The value $a$ is called the residue of $\lambda$ around $\Sigma$. 
\end{definition}

The form $\lambda_0$ defined in section \ref{s:sdbdef} on $\sdb(Z(s_\eps),\nf 1k)^*$ provides a Liouville form in $\op(Z(s_\eps))^*$ that is tame along $Z(s_\eps)$, with residue $-\frac 1k$ (recall our notation $\op(X)^*=\op(X)\priv X$).  On the other hand, the form $\lambda_\eps$, defined in the previous section on the larger set $\sdb(N,\tau)\priv Z(s_\eps)$, may not be tame. 
In the K\"ahler setting,  $\lambda_\eps$ is  $-\frac 1k d^c\ln|s_\eps|$, which is tame, but the proof I am aware of does not  pass easily to the quasi-holomorphic setting. The next two lemma allow to modify $\lambda_\eps$ locally near $Z(s_\eps)$ to a tame Liouville form, keeping its skeleton fixed.  
We recall that  $s_\eps$ vanishes transversally, not on $N_\infty$, so $Z(s_\eps)\subset \sdb(N,\tau)$ is a smooth symplectic submanifold. 
As a  zero-set of a section of $\cl^k\to \spb(N,\tau)$, it is Poincaré-dual to $k[\om_0]$, so it has a \nbd that can be identified with a \nbd of the zero-section of a symplectic disc bundle, on which we have a radial coordinate $P$ and a Liouville form $\lambda_0^\eps$ whose vector field is negatively proportional to $\frac\partial{\partial P}$. 
\begin{lemma}\label{l:lambdagrowsR}
For $P$ small enough, $d |s_\eps|^2\big(\frac\partial{\partial P}\big)>0$. 
\end{lemma}
\noindent{\it Proof:} Since $s_\eps$ vanishes transversally in the case at hand, $|s_\eps|^2$ is Morse-Bott along $Z(s_\eps)$ with positive definite  transverse Hessian. Thus, in local coordinates $(z,P,\theta)$ provided by the symplectic disc bundle model, and denoting $P=\rho^2$, we have 
$$
\frac{\partial^2 |s_\eps|^2}{\partial \rho^2}(z,0,\theta)>0 \hspace{1cm} \forall z,\forall \theta,
$$
so this second derivative is positive also for sufficiently small $\rho$. Since  the first derivative in the radial direction has to vanish at $\rho=0$ because $Z(s_\eps)$ is a minimum of the function $|s_\eps|^2$, this first radial derivative is positive on a small enough pointed \nbd of $Z(s_\eps)$. 
\cqfd

\begin{lemma}\label{l:notametotame} Let $(M,\om)$ be a possibly open symplectic manifold, $\Sigma$ a smooth symplectic divisor locally Poincaré-dual to $\ell[\om]$, $\pi:\op(\Sigma)\to \Sigma$ the natural projection when $\op(\Sigma)$ is identified with a \nbd of the zero section in $\sdb(\Sigma,\om_{|\Sigma},\ell)$, and $f:\op(\Sigma)\to [0,1]$ a smooth function whose zero-set is $\Sigma$. Let $\lambda,\lambda'$ be two Liouville forms defined in $\op(\Sigma)^*$ that satisfy $X_\lambda\cdot f>0$ and $X_{\lambda'}\cdot f>0$. Then for all $0<\eps<\eps'\ll 1$, there exist a closed form $\vartheta$ on $\Sigma$ and a Liouville vector field $\lambda''$ in $\op(\Sigma)^*$ that coincides with $\lambda+\pi^*\vartheta$ in $\{f<\eps\}$, with $\lambda'$ in $\{f>\eps'\}$, and such that $X_{\lambda''}\cdot f>0$.   
\end{lemma}
\noindent{\it Proof:}  $\lambda'-\lambda$ is a closed $1$-form on $\op(\Sigma)^*$, that we see as a pointed \nbd of a symplectic disc bundle. Let $\gamma(\theta):=(z_0,R_0,\theta)\in \op(\Sigma)^*$ be a small loop in a fiber. Since $H_1(\op(\Sigma)^*)=H_1(\Sigma)+\la \gamma \ra$, we can write 
$$
\lambda'-\lambda=\pi^*\vartheta+\nu,
$$
where $\vartheta$ is a smooth closed $1$-form on $\Sigma$ and $\nu$ is a closed $1$-form that vanishes on $H_1(\Sigma)$. Consider a $2$-cycle $S\subset \Sigma$,  deformed to $S'\subset \op(\Sigma)$ that intersects $\Sigma$ transversally. By Stokes theorem,  
$$
-c_1(\cn \Sigma,[S])[\lambda'-\lambda](\gamma)=\int_{S'} \om-\om=0,
$$
so $\nu$ is in fact exact. Thus we can write 
$$
\lambda'-\lambda=\pi^*\vartheta+dh,
$$
for some smooth function $h:\op(\Sigma)^*\to \R$. Let now $\eps<\eps'$ be such that $\{f<\eps'\}\subset \op(\Sigma)$, $\chi:\R\to [0,1]$ be a smooth function that equals $0$ before $\eps$ and $1$ after $\eps'$. Define 
$$
\lambda'':=\lambda+\pi^*\vartheta+d(\chi (f) h)=\chi(f)\lambda'+(1-\chi(f))\lambda+hd\chi\circ f. 
$$
The first formula shows that $\lambda''$ is again a Liouville form that coincides with $\lambda+\pi^*\vartheta$ on $\{f<\eps\}$ and with $\lambda'$ on $\{f>\eps'\}$. The second shows that the associated Liouville vector field is 
a positive barycenter of $X_\lambda$ and $X_\lambda'$ modulo a vector field $Z$ which is tangent to the level sets of $f$. This implies that 
$X_{\lambda''}\cdot f>0$, as claimed. \cqfd

Let us now come back to our framework: $\lambda_\eps$ is a Liouville form in the complement of the smooth divisor $Z(s_\eps)\subset \sdb(N,\tau)$, this divisor being Poincaré-dual to $k[\om_0]$ in $M=\spb(N,\tau)$, hence locally around itself. 
We have a radial coordinate $P$ in $\op(Z(s_\eps))$ and a local tame Liouville form $\lambda_0^\eps$, whose associated Liouville vector field is $(P-\frac 1k)\frac\partial{\partial P}$, coming from the identification of $\op(Z(s_\eps))$ with a symplectic disc bundle. By lemma \ref{l:lambdagrowsR}, $d|s_\eps|^2\big(\frac\partial{\partial P}\big)>0$ so  $X_{\lambda_0^\eps}\cdot (-|s_\eps|^2)>0$ in $\op(Z(s_\eps))^*$ and the same positivity holds for $X_\eps$ by the quasi-holomorphicity of $s_\eps$ (lemma \ref{l:qhsec}). By lemma \ref{l:notametotame}, there exists a Liouville form $\lambda_\eps''$ on $\sdb(N,\tau)\priv Z(s_\eps)$ that coincides with $\lambda_\eps$ outside a \nbd of $Z(s_\eps)$ and with $\lambda_0^\eps+\pi^*\vartheta$ on $\op(Z(s_\eps))^*$, where $\vartheta$ is a smooth $1$-form on $\op(Z(s_\eps))$. The latter form is tame along $Z(s_\eps)$, so $\lambda_\eps''$ is tame. Moreover, this same lemma guarantees that its associated Liouville vector field $X_\eps''$ is still gradient-like for $-\ln|s_\eps|$. Since moreover $Z(s_\eps)$ does not meet $N_\infty\cup \skel(s_\eps)$,  $X_\eps''$ points outside $N_\infty$ just the same as $\lambda_\eps$ and has the same skeleton.

In conclusion, we have obtained the following (we drop the comas and the subscript from $\lambda''_\eps$):

\begin{proposition}\label{p:skels} 
Let $(\cp,\alpha)\to (N,\tau)$ be  a prequantization bundle. Let $(\cl^k,\nabla)\to (N,\tau)$ be a Hermitian bundle endowed with a connection of curvature $-ik\tau$, $\sigma_k:N\to \cl^k$ a quasi-holomorphic section that vanishes transversally, and $\Gamma_k$ its associated skeleton:
$$
\Gamma_k:=\bigcup_{z\in \crit \ln |\sigma_k|} W^s(X_{\sigma_k},z). 
$$ 
Let $\Lambda_k\subset \cp$ be a $k$-fold isotropic lift of $\Gamma_k$. Then there exists a smooth symplectic divisor $Z_k\subset (\sdb(N,\tau),\om_0)$ and a Liouville form $\lambda_k$ on $\sdb(N,\tau)\priv Z_k$, tame along $Z_k$  with residue $-\frac 1k$, whose Liouville vector field $X_k$ parts $\sdb(N,\tau)$ into two subsets:
\begin{itemize}
\its the points whose $X_k$-trajectories abute in finite time to $Z_k$. This subset, denoted $\cb(Z_k,X_k)$, is called the basin of attraction of $Z_k$ under $X_k$. It is symplectomorphic to $\sdb(Z_k,\om_{0|Z_k},k)$. 
\its the points that never reach $Z_k$, on which the flow of $X_k$ is defined over $\R$.  This subset, called the skeleton of $X_k$ and denoted $\skel(Z_k,\lambda_k)$ is explicit: 
$$
\skel(Z_k,\lambda_k)=\overline{\Phi^{-1}(\Lambda_k\times (0,1))},
$$ 
$\Phi$ being the identification between $\sdb(N,\tau)\priv N_0$ with $S^{< 0}\cp$ described in lemma \ref{l:sdbsp}.  
\end{itemize}
\end{proposition}
\noindent{\it Proof:} Two claims need a justification beyond the discussion above. The first one is the fact that 
$$
\cb(Z_k,X_k)\overset\om\simeq \sdb(Z_k,\om_{0|Z_k},k) 
$$
and holds by \cite{opshtein2} because $\lambda_k$ is tame along $Z_k$ with residue $-\nf 1k$. 
The second is the fact that  
$\Phi^{-1}(\Lambda_k\times(0,1))$ is the skeleton of a Liouville vector field {\it for any} $k$-fold isotropic lift $\Lambda_k$ of $\Gamma_k$: so far we have only obtained {\it one specific} $\lambda_{\text{spec}}$ whose skeleton is {\it one such specific} subset, that we denote $\Phi^{-1}(\Lambda_{\text{spec}}\times (0,1))$. But the uniqueness part of lemma \ref{l:skelift} shows that all these $k$-fold lifts differ by a global rotation of $\sdb(N,\tau)$, so any such $\Lambda_k$ can be written $\cR_\alpha(\Lambda_{\text{spec}})$ for some angle $\alpha$. Since $\cR_\alpha$ is a global symplectomorphism of $\sdb(N,\tau)$, $\lambda:=\cR_\alpha^*\lambda_{\text{spec}}$ is a Liouville form tame along $Z:=\cR_\alpha(Z(s_\eps))$ with the same residue $-\frac 1k$ as $\lambda_{\text{spec}}$, whose skeleton is $\cR_\alpha(\Phi^{-1}(\Lambda_{\text{spec}}\times (0,1)))=\Phi^{-1}(\Lambda_k\times (0,1))$. \cqfd

\section{Proof of theorem \ref{t:legbarp}.} 
Let $(N,\tau)$ be a closed symplectic manifold with integral symplectic class and $(\cp,\alpha):=\cp(N,\tau)$ its prequantization bundle. Let $\Gamma_k\subset N$ be an isotropic skeleton associated to a quasi-holomorphic section of $\cl^k\to (N,\tau)$. Let $\Lambda_k$ be a $k$-fold Legendrian lift of $\Gamma_k$ to $(\cp,\alpha)$ described in corollary \ref{c:legbardot1}. Let also $Z_k$ the symplectic divisor of $(\sdb(N,\tau),\om_0)$ and $\lambda_k$ the Liouville form on its complement given by proposition \ref{p:skels}, with the associated decomposition 
$$
\begin{array}{l}
\sdb(N,\tau)=\cb(Z_k,\lambda_k)\sqcup \skel(Z_k,\lambda_k),\\
\skel(Z_k,\lambda_k)=\overline{\Phi^{-1}(\Lambda_k\times(0,1))},\\
\cb(Z_k,\lambda_k)\overset \om\approx\sdb(Z_k,\om_{0|Z_k},k). 
\end{array}
$$

Let  $H_t:M\to [1,+\infty)$ be a contact Hamiltonian and $\Lambda$ any closed Legendrian submanifold of $(\cp,\alpha)$. Assume that there is no $H$-chord between $\Lambda$ and $\Lambda_k$ of length $\leq T$ and that 
$$
\bigcup_{t\in [0,T]}\phi_H^t(\Lambda)
$$
 is  embedded. Then by proposition \ref{p:laglift} there exists a closed Lagrangian submanifold $L\approx \Lambda\times S^1$ in $S^{<0}\cp\priv (\Lambda_k\times (0,1))$ with $\ca_{\min}(L)= T$. By lemma \ref{l:sdbsp}, $L$ can be seen as a Lagrangian submanifold of $\sdb(N,\tau)\priv N_0$ in the complement of 
$\Phi^{-1}(\Lambda_k\times (0,1))$, hence in $\cb(Z_k,\lambda_k)\approx \sdb(Z_k,\om_{0|Z_k},k)$  
with the same minimal area $T$. Since $L$ is closed, it lies in a compact subset of $\sdb(Z_k,\om_{0|Z_k},k)$, 
 which compactifies to a closed symplectic uniruled manifold  of width $<\frac 1k$ by Biran's compactification procedure described in the begining of \S \ref{s:spb} (see also \cite{biran2}). The latter is a symplectic sphere bundle $(Q,\om)$ over $Z_k$, whose fibers have area $<\frac 1k$, possesses two disjoint symplectic sections $Z_k^0$ and $Z_k^\infty$. The Lagrangian submanifold $L$ lies in the complement of $Z_k^\infty$, and $Z_k^0$ is symplectomorphic to $Z_k$, which was Poincaré-Dual to $k[\om_0]$. Since $[\om_0]$ is an integral class, $[\om]$ takes values in $k\Z$ on $H_2(Z_k^0,\R)$, so those of these values that are positive all exceed the area of the fiber. By theorem \ref{t:cimo}, there exists a symplectic disc in 
 $$Q\priv Z_k^\infty\subset \sdb(Z_k,\om_{0|Z_k},k)\subset \sdb(N,\tau)$$
 with boundary in $L$ and symplectic area $<\nf 1k$. Thus,
$$
\ca_{\min}(L)=T< \frac 1k. \hspace*{2cm}\square
$$


{\footnotesize
\bibliographystyle{alpha}
\bibliography{biblio}
}

\end{document}